\documentclass[11]{article}

\usepackage{amsmath,amsfonts,amsthm}

\usepackage{graphicx,color,pict2e}

\definecolor{refkey}{rgb}{1, 0.5, 0}  
\definecolor{labelkey}{rgb}{1,1,1}
\definecolor{labelkey}{rgb}{0.1,1,0.2}

\newtheorem{remark}{Remark}[section]

\newcommand{\abs}[1]{\left| #1 \right|}

\def\bel{\begin{equation}\label}
\def\eeq{\end{equation}}
\def\bega{\begin{array}}
\def\enda{\end{array}}

\def\ve{\varepsilon}

\title{Global Riemann Solvers for Several $3\times3$ Systems of Conservation Laws with Degeneracies}
\author{Wen Shen\\
{Mathematics Department, Pennsylvania State University,}\\ 
{University Park, PA 16802, USA.}}

\begin{document}

\maketitle

\begin{abstract} 
We study several $3\times 3$ systems of conservation laws,
arising in modeling of two phase flow with rough porous media
and traffic flow with rough road condition. 
These systems share several features. 
The systems are of mixed type, with various degeneracies.
%along some curves in the domain  
%the systems are parabolic degenerate,
%where eigenvalues and eigenvectors of different families coincide. 
%waves from different families switch the ordering of their speed. 
%This causes nonlinear resonance. 
Some families are linearly degenerate, while others are not genuinely nonlinear.
Furthermore, along certain curves in the domain  
the eigenvalues and eigenvectors of different families coincide. 
Most interestingly, in some suitable Lagrangian coordinate
the systems are partially decoupled, where some unknowns can be solved
independently of the others.
Finally, in special cases the systems reduce to 
some $2\times2$ models, 
which have been studied in the literature. 
Utilizing  the insights gained from these features, 
 we construct global Riemann solvers for all these models. 
Possible treatments on the Cauchy problems are also discussed. 
\end{abstract}

%%%%%%%%%%%%%%%%%%%%%%%%%%%
\section{Introduction} \label{sec0}
%%%%%%%%%%%%%%%%%%%%%%%%%%%

Scalar conservation laws with discontinuous flux functions have attracted significant 
research interests in recent years,  and exciting progresses have been made. 
See for example a survey paper \cite{And} and references therein. 
In a general setting, a scalar conservation law
\begin{equation}\label{a.1}
 u_t + g(a(x),u)_x=0
 \end{equation}
where $a(x)$ contains discontinuity, can be written into a $2\times 2$ system,
by adding a trivial equation for $a(x)$:
\begin{equation}\label{a.2}
 \begin{cases}
  u_t + g(a,u)_x =0,\\
  a_t=0.
 \end{cases}
\end{equation}

For the general triangular system~\eqref{a.2}, when $g_u(a,u)=0$, 
the two eigenvalues and eigenvectors of the two families coincide,
and the system is not hyperbolic. 
In the literature this is referred to  as   parabolic degeneracy. 
Utilizing the vanishing viscosity solution of 
\begin{equation}\label{a.5}
 \begin{cases}
  u_t + g(a,u)_x =\ve u_{xx},\\
  a_t=0,
 \end{cases}
\end{equation}
as $\ve\to 0+$, 
solutions of Riemann problems can be uniquely determined.
Such admissible condition for jumps in $a(x)$ 
leads to the {\em minimum jump condition}. 
See~\cite{GimseRisebro} and  some more recent works~\cite{GS2017, ShenN}.

Triangular systems~\eqref{a.2} arises in many physical models.
Take for example  the two phase flow models in reservoir simulations.
Consider a simple polymer flooding model with single component
\begin{equation}\label{a.3}
 \begin{cases}
  s_t + f(s,c)_x =0,\\
  (cs)_t + (c f(s,c))_x=0.
 \end{cases}
\end{equation}
Here, $s\in[0,1]$ is the saturation of the water phase,  
$c\in[0,1]$ is the fraction of polymer dissolved in water, 
and $f(s,c)$ is the fractional flow for the water phase. 
One assumes uniform porous media, no gravitation force, and no adsorption 
of the polymer into the porous media. 
Introducing a Lagrangian coordinate $(\tau,y)$ (see~\cite{Wagner}), with 
$$ y_x =s, \quad y_t =-f, \quad y(0,0)=0, \qquad \tau=t,$$
the system~\eqref{a.3} can be written as a triangular 
system
\begin{equation}\label{a.4}
 \begin{cases}
  (1/s)_\tau - (f(s,c)/s)_y =0 , \\
  c_\tau=0.
 \end{cases}
\end{equation}

In this paper, 
we consider the two-phase polymer flooding in rough media,
where the permeability function of 
the porous media may be  discontinuous. 
Let $k(x)$ be the absolute permeability of the rock, 
system~\eqref{a.3} is extended to the following  
$3\times 3$ systems of conservation laws, 
where we also take into consideration of the adsorption
of polymers into the rock:
\begin{align}
 s_t + f(s,c,k)_x ~=~0 ,\label{0.1} \\
 (m(c)+cs)_t + (cf(s,c,k))_x ~=~0,\label{0.2} \\
 k_t ~=~0.\label{0.3} 
\end{align}
Here, the unknown vector is $(s,c,k)$, and 
the function $m(c)$ describes the adsorption of polymer in the porous media. 
We assume that $m$ depends only on $c$.

The main objective of this paper is the construction of  
global Riemann solvers for~\eqref{0.1}-\eqref{0.3}, 
under 3 different situations:
\begin{itemize}
\item  We neglect the gravity force and the adsorption. See section~\ref{sec2};
\item We consider the adsorption and neglect the gravity force. See section~\ref{sec3};
\item We consider the gravity force and neglect the adsorption. See section~\ref{sec5}.
\end{itemize}

As an additional model, we also treat a 
$3\times 3$ system modeling traffic flow in section~\ref{sec4}.
The traffic flow system has some similar features to the polymer flooding models. 

We remark that, global Riemann solvers for general 
nonlinear systems of hyperbolic conservation laws can not always be constructed,
due to the nonlinearity of the flux function. 
Such a task is possible here, thanks to the special properties of the models.
Once a Riemann solver is available, remarks are given on possible approaches
to establishing existence of solutions for Cauchy problem, 
for some of the cases. 
Finally, concluding remarks are given in section~\ref{sec6}
where more future works are suggested.

%%%%%%%%%%%%%%%%%%%%%%%%%%%
\section{A simple model for polymer flooding in two phase flow with rough media}\label{sec2}
\setcounter{equation}{0}
%%%%%%%%%%%%%%%%%%%%%%%%%%%%%

We first consider the two phase flow mode of polymer flooding~\eqref{0.1}-\eqref{0.3},
where we neglect the adsorption effect and the gravitation effect, i.e., 
\begin{align}
 &s_t + f(s,c,k)_x =0, \label{1.1} \\
& (cs)_t + (cf(s,c,k))_x =0,\label{1.2} \\
 &k_t =0.\label{1.3} 
\end{align}
The flux function $f(s,c,k)$  has the following properties.  
For any given $(c,k)$, the mapping $s\mapsto f$ is the famous S-shaped 
Buckley-Leverett function~\cite{BL}
with a single inflection point. 
We have 
$$f(s,c,k)\in[0,1], \qquad f_s(s,c,k)\ge0, \qquad \mbox{for all } (s,c,k),$$ 
and 
\begin{equation}\label{fconds}
 f(0,c,k)=0, \qquad f(1,c,k)=1, \qquad f_s(0,c,k)=0, \qquad 
f_s(1,c,k) =0,\qquad \forall(c,k).
\end{equation}
Furthermore, it's physically  reasonable to assume that the flux 
decreases with more dissolved  polymer, and increases with increasing permeability,
i.e.,
\begin{equation}\label{fconds2}
  f_c(s,c,k) < 0, \qquad f_k(s,c,k) >0, \qquad \forall (s,c,k).
\end{equation}
The assumptions~\eqref{fconds2} 
simplify the analysis, allowing clearer presentation of the main ideas. 
We remark  that,
if we remove the assumptions~\eqref{fconds2}, 
a similar analysis can be carried out, but with heavier details. 

%%%%%%%%%%%%%%%%%%%%%%%%%%
\subsection{Riemann solver for the reduced $2\times2$ model}\label{sec2.1}
%%%%%%%%%%%%%%%
Observe that when $k$ is constant, the system~\eqref{1.1}-\eqref{1.3}
reduces to a $2\times 2$ system
\begin{align}
 &s_t + f(s,c)_x =0, \label{1.10} \\
& (cs)_t + (cf(s,c))_x =0.\label{1.20} 
\end{align}
With a slight abuse of notation we write $f(s,c)=f(s,c,k)$ when $k$ is a constant. 
The Jacobian matrix of the flux function is triangular
$$ 
J = \begin{pmatrix} f_s(s,c) & f_c(s,c) \\
0 & f(s,c)/s  
\end{pmatrix}.
$$
The two eigenvalues and the corresponding right eigenvectors of $J$ are 
$$\lambda_s= f_s(s,c), \quad \lambda_c= f(s,c)/s, \qquad 
r_s=\begin{pmatrix} 1 \\ 0 \end{pmatrix}, \quad
r_c=\begin{pmatrix} f_c(s,c) \\ \lambda_s-\lambda_c \end{pmatrix}.
$$
When $\lambda_s=\lambda_c$, the two eigenvectors also coincide,
therefore the system becomes parabolic degenerate. 
Since the difference $\lambda_s-\lambda_c$ can change sign, 
nonlinear resonance occurs, and the total variation of the unknown
can blow up in finite time, see~\cite{Temple, Temple2}.  
Therefore weak solutions $(s,c)$ are not 
defined in the class of BV functions. 

System~\eqref{1.10}-\eqref{1.20}  has been studied in quite some detail in the literature.
It is known that 
Riemann problems  for~\eqref{1.10}-\eqref{1.20} can be solved globally, 
generating entropy solutions that 
are the vanishing viscosity limit, see~\cite{GS2017}. 
See also a recent work~\cite{WS}, where Riemann problems, as well as 
the existence of solutions for the Cauchy problems are treated, with the consideration
of the gravity force. 

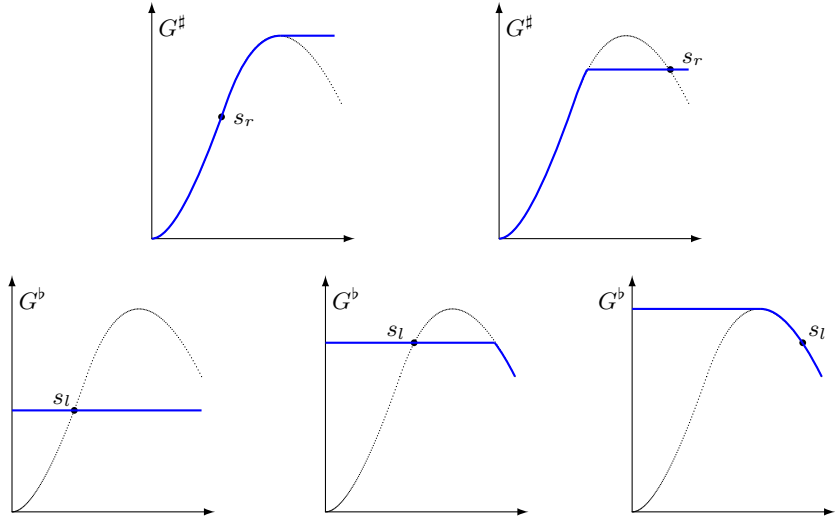
\begin{figure}[htbp]
\begin{center}
\setlength{\unitlength}{0.9mm}
% upper left ------------
\begin{picture}(50,38)(0,0)  
\put(0,0){\vector(1,0){30}}
\put(0,0){\vector(0,1){35}} \put(1,30){\small$G^\sharp$}
\qbezier[80](0,0)(4,0)(11,20)
\qbezier[80](11,20)(18,40)(28,20)
\put(10.3,18){\circle*{1}}\put(12,17){\small $s_r$}
\color{blue}
\thicklines
\qbezier(0,0)(4,0)(11,20)
\qbezier(11,20)(14.5,30)(19,30)
\put(19,30){\line(1,0){8}}
\end{picture}
\begin{picture}(30,38)(0,0)  % upper right ----------
\put(0,0){\vector(1,0){30}}
\put(0,0){\vector(0,1){35}}
\put(1,30){\small$G^\sharp$}
\qbezier[80](0,0)(4,0)(11,20)
\qbezier[80](11,20)(18,40)(28,20)
\thicklines
\put(25.3,25){\circle*{1}}\put(27,26){\small $s_r$}
\color{blue}
\qbezier(0,0)(4,0)(11,20)
\qbezier(11,20)(12,23)(13,25)
\put(28,25){\line(-1,0){15}}
\end{picture}

%
% lower left -------
\begin{picture}(45,40)(0,0)  
\put(0,0){\vector(1,0){30}}
\put(0,0){\vector(0,1){35}}
\put(1,30){\small$G^\flat$}
\qbezier[80](0,0)(4,0)(11,20)
\qbezier[80](11,20)(18,40)(28,20)
\thicklines
\put(9.2,15){\circle*{1}}\put(6,16){\small $s_l$}
\color{blue}
\put(28,15){\line(-1,0){28}}
\end{picture}
% lower middle ----
\begin{picture}(40,38)(0,0)  
\put(0,0){\vector(1,0){30}}
\put(0,0){\vector(0,1){35}}
\put(1,30){\small$G^\flat$}
\qbezier[80](0,0)(4,0)(11,20)
\qbezier[80](11,20)(18,40)(28,20)
\thicklines
\put(13.1,25){\circle*{1}}\put(9,26){\small $s_l$}
\color{blue}
\qbezier(25,25)(26.5,23)(28,20)
\put(0,25){\line(1,0){25}}
\end{picture}
% lower right -------------
\begin{picture}(35,38)(-4,0)  
\put(0,0){\vector(1,0){30}}
\put(0,0){\vector(0,1){35}}
\put(-5,30){\small$G^\flat$}
\qbezier[80](0,0)(4,0)(11,20)
\qbezier[80](11,20)(18,40)(28,20)
\thicklines
\put(25.2,25){\circle*{1}}\put(26,26){\small $s_l$}
\color{blue}
\qbezier(19,30)(23,30)(28,20)
\put(0,30){\line(1,0){19}}
\end{picture}
\caption{The functions $G^\sharp(s;s_r)$ and $G^\flat(s;s_l)$ for various cases of $s_r$ and $s_l$. The dotted curve is the graph of $s\mapsto g(s,c)$ for a fixed $c$.  }
\label{fig:G}
\end{center}
\end{figure}

We briefly summarize the Riemann solver for~\eqref{1.10}-\eqref{1.20},
which will be useful for the solution of the full $3\times3$ system. 
Given the Riemann data $(s_l,c_l), (s_r,c_r)$, 
we define the functions
$$g(s,c) ~\dot=~ f(s,c)/s,\qquad
g_l(s) ~\dot=~ g(s,c_l), \qquad g_r(s) ~\dot=~ g(s,c_r),
$$
and the monotone functions
\begin{align}
G^\sharp(s;s_r) &~\dot=~ \begin{cases}
\max\left\{ g_r(s); s\in[s_r,s]\right\}, & \mbox{if}~s\ge s_r,\\
\min\left\{ g_r(s); s\in[s,s_r]\right\}, & \mbox{if}~s\le s_r,
\end{cases}\label{Gsharp}
\\[1mm]
G^\flat(s;s_l) &~\dot=~ \begin{cases}
\max\left\{ g_l(s); s\in[s_l,s]\right\}, & \mbox{if}~s\ge s_l,\\
\min\left\{ g_l(s); s\in[s,s_l]\right\}, & \mbox{if}~s\le s_l.
\end{cases} \label{Gflat}
\end{align}
See plots in Figure~\ref{fig:G} for illustrations.

Note that the mapping $s\mapsto G^\sharp$ is  increasing, 
while $s\mapsto G^\flat$ is  decreasing.
For any given $(s_l,s_r)$, 
there exists a unique $\hat G$ value where the graphs of the two mappings
cross each other. 
Let $\sigma$ denote the speed of the $c$~jump, and 
let $s_\pm$ denote the traces $s(t,\sigma t\pm)$ in the Riemann solution.
Then, $s_\pm$ are determined as the {\em minimum jump path that connects
$G^\sharp$ and $G^\flat$}, (see also~\cite{GimseRisebro}),
where we have 
\begin{equation}\label{gGGg}
g_l(s_-)=G^\flat(s_-;s_l)=G^\sharp(s_+;s_r)=g_r(s_+) ~\dot=~\sigma.
\end{equation}
Then, the solution of the Riemann problem is obtained by patching together 
the solution of 
\begin{equation}\label{Rl}
 s_t + f(s,c_l)=0, \qquad
s(0,x)=\begin{cases} s_l & (x<0),\\   s_- & (x>0), \end{cases}
\end{equation}
for $x<\sigma t$, and  for $x> \sigma t$ the  solution of
\begin{equation}\label{Rr}
 s_t + f(s,c_r)=0, \qquad
s(0,x)=\begin{cases} s_+ & (x<0),\\   s_r & (x>0) .\end{cases}\end{equation}

%%%%%%%%%%%%%%%%%%%%%%%%%%%
\subsection{The Lagrangian coordinate}
%%%%%%%%%%%%%%%%%%%%%%%%%

Define the Lagrangian coordinate $(\phi,\psi)$  (introduced in \cite{Splitting})
\begin{equation}\label{eq:c}
 \phi_x = -s, \quad \phi_t = f(s,c,k),  \quad \phi(0,0)=0, \qquad
\psi=x.
\end{equation}
Here one can interpret $\phi$ as the potential for the first equation. 
In fact, for any $(t,x)$,  the value  $\phi(t,x)$ denotes the line integral
$$ \phi(t,x) = \int_{(0,0)}^{(t,x)}  f (s,c,k)\, dt - s \, dx.$$
Thanks to \eqref{1.1}, this line integral is path independent. 
Assuming $s>0$ so that $f>0$,  the coordinate change is well defined.
In this Lagrangian coordinate,  the system \eqref{1.1}-\eqref{1.3} takes the form 
\begin{align}
 \left(\frac{1}{f}\right)_\psi-\left(\frac{s}{f}\right)_\phi =0, \label{1.6}\\
c_\psi =0,\label{1.7}\\
k_\phi=0.\label{1.8}
\end{align}
Note that the second and third equations are decoupled. 
Since $k$ is a material parameter, the decoupling is not surprising.
The decoupling for $c$ indicates that the thermo-dynamics (governed by~\eqref{1.7})
is independent of the hydro-dynamics (governed by~\eqref{1.6}).
This is the most interesting feature of the model. 
It implies that, in the $(\phi,\psi)$ coordinates,  
$k$ is constant along lines parallel to $\phi$-axis,
and $c$ is constant along lines parallel to the $\psi$-axis. 

\begin{figure}[htbp]
\begin{center}
\setlength{\unitlength}{1.1mm}
\begin{picture}(80,42)(-25,-16)  
\put(-20,0){\vector(1,0){60}}\put(41,0){$\psi$}
\put(0,-15){\vector(0,1){35}}\put(0,21){$\phi$}
\thicklines
\color{blue}
\put(0,0){\line(-1,2){10}} 
\put(0,0){\line(1,-1){18}}
\put(-25,-5){initial condition}\put(-18,-8){ at $t=0$}
\thinlines
\put(-5,-2){\vector(0,1){10}}
\put(-2,-7){\vector(1,0){8}}
\color{cyan}
\put(10,-10){\vector(1,1){5}}\put(14,-10){$t$  increases}
\color{red}
\put(-3,6){\line(1,0){30}}\put(30,5){\small $c$=constant}
\put(3,-3){\line(1,0){30}}\put(35,-4){\small $c$=constant}
\color{green}
\put(4,-4){\line(0,1){20}} \put(4,17){\small $k$=constant}
\put(8,-8){\line(0,1){24}} 
\end{picture}
\caption{The connection between the two coordinate for Riemann solutions.}
\label{fig:1}
\end{center}
\end{figure}
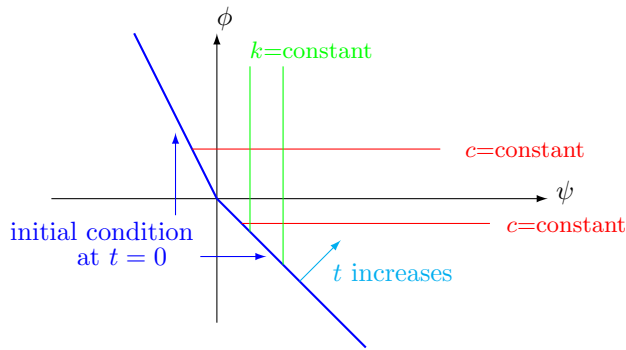

We illustrate the coordinate change in Figure~\ref{fig:1}, with  Riemann data $s_l,s_r>0$. 
The line $t=0$  now consists of two rays from the origin, 
in the Lagrangian coordinate $(\phi,\psi)$, indicated in blue in Figure~\ref{fig:1}.
For general initial data, the line for $t=0$ will be replaced by a curve, continuous and decreasing,
but might not be differentiable everywhere. 
Here we see clearly how the values of $k$ and $c$ are ``brought into''
the region $t>0$ from the initial condition.

Note that when $s=0$ then $f=0$, and the conserved quantity for the equation 
in Lagrangian coordinate blows up to infinity. This is when we have ``vacuum''. 
If $s(0,x)=0$ for an intervals $[x_1,x_2]$, then the blue curve 
in Figure~\ref{fig:1} will have a horizontal line segment. 
Since $c$ has no meaning when $s=0$, we may assign the $c$ value
along this segment as a linear function
connecting $c_1=c(0,x_1)$ and $c_2=c(0,x_2)$. 
If $c_1\not=c_2$, then $c$ has a jump in its solution.

This illustration  shows that, once the initial data is given, 
the values for $(k,c)$ are known at every point $(\psi,\phi)$. 
In particular, if $(c,k)$ are  initially smooth and $s(0,x)>0$, 
then $(c,k)$ remain smooth forever; 
if they contains discontinuities initially, then they will be discontinuous for $t>0$. 

We remark that this Lagrangian coordinate is different from the one used by 
Wagner in his seminal paper~\cite{Wagner}, 
for Euler's equation.
If we apply Wagner's Lagrangian coordinate to~\eqref{1.1}-\eqref{1.3},
 only the equation for $c$ will be decoupled. This doesn't offer 
 the same insight.

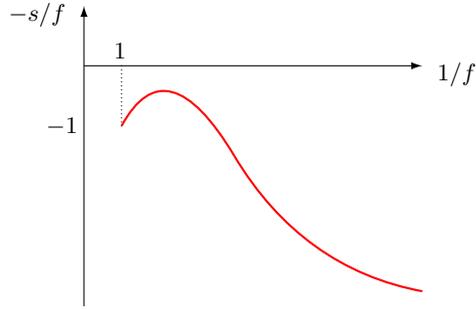
\begin{figure}[htbp]
\begin{center}
\setlength{\unitlength}{1mm}
\begin{picture}(40,44)(0,0)  
\put(-5,0){\vector(0,1){40}}\put(-15,38){\small $-s/f$}
\put(-1,33){\small 1}
\put(-10,23){\small $-1$}
\put(-5,32){\vector(1,0){45}}\put(42,30){\small $1/f$}
\qbezier[15](0,24)(0,28)(0,32)
\thicklines
\color{red}
\qbezier(0,24)(6,35)(15,20)\qbezier(15,20)(24,5)(40,2)
\end{picture}
\caption{The graph for the mapping $(1/f)\mapsto (-s/f)$.}
\label{fig:fs}
\end{center}
\end{figure}

We  now consider  a scalar conservation law with possibly discontinuous
coefficients
\begin{equation}\label{eqL}
 \left(\frac{1}{f(s;c,k)}\right)_\psi- \left(\frac{s}{f(s;c,k)}\right)_\phi=0, 
 \end{equation}
where $c,k$ are given functions, possibly discontinuous.
A typical plot of the ``flux'' function, i.e., the graph for 
the mapping $(1/f )\mapsto(- s/f)$ is shown in Figure~\ref{fig:fs}.

We remark that scalar conservation laws with horizontal and vertical discontinuities
was studied in~\cite{CocliteRisebro}, where Riemann problem and Cauchy problem are 
studied, under suitable assumption of the flux function. 
We speculate that an extension of~\cite{CocliteRisebro}  could provide 
existence and well-posedness for the Lagrangian system~\eqref{1.6}-\eqref{1.8}. 
Details may be worked out in future works. 
Furthermore, it would also be interesting to obtain equivalent results directly 
for the Eulerian system,  see Remark~\ref{rm1}.

%%%%%%%%%%%%%%%%%%%%%%
\subsection{Wave properties and a global Riemann solver} 
\label{sec2.3}
%%%%%%%%%%%%%%%%%%%%%%

We now return to the full Eulerian system~\eqref{1.1}-\eqref{1.3}. 
Treating $(s,c,k)$ as the unknown vector, 
the Jacobian matrix of the flux function is triangular: 
$$
J=\begin{pmatrix} f_s & f_c & f_k \\
0 & f/s & 0 \\
0 & 0 &0
\end{pmatrix}.
$$
Naming the three families as the $s$, $c$ and $k$ families, 
we have the following 3 eigenvalues 
$$
\lambda_s= f_s, \qquad 
\lambda_c=f/s,\qquad
\lambda_k=0,
$$
and three corresponding right eigenvectors, 
$$
r_s=\begin{pmatrix}1 \\ 0 \\ 0 \end{pmatrix},
\qquad
r_c= 
\begin{pmatrix} -f_c \\ f_s-f/s \\  0  \end{pmatrix}, 
\qquad 
r_k=\begin{pmatrix} -f_k\\ 0 \\ f_s\end{pmatrix}.
$$

Direct computations give: 
$$ 
\nabla \lambda_s \cdot r_s = f_{ss} , \qquad
\nabla \lambda_c\cdot r_c =0, \qquad \nabla \lambda_k \cdot r_k =0.
$$
Thus, the $c$ and $k$ families are linearly degenerate,
where discontinuities are all contacts, and shock curves coincide with 
rarefaction curves. 
For the $s$ family, 
since $f_{ss}$ changes sign, the family is not genuinely nonlinear.
However, the integral curves for $s$ family are straight lines where $c,k$
are both constants.
These integral curves also coincide with the $s$~shock curves,
make it easier to find waves of the $s$~family.
Note also that, 
when $f_s=f/s$, we have $\lambda_s=\lambda_c$ and $r_c=r_c$,
and the system is parabolic degenerate, where nonlinear resonance occurs.
In summary,  the system~\eqref{1.1}-\eqref{1.3} is of Temple class,
but of mixed type with degeneracies. 

By the Rankine-Hugoniot jump conditions, we have the following wave properties:
\begin{itemize}
\item The $k$~wave is the slowest which travels with speed 0.  Along any $k$ wave, 
the functions $f, c$  are continuous;
\item The $c$~wave travels with positive speed. Crossing it, $f/s, k$ remain continuous;
\item The $s$~wave  travels with positive speed. Crossing it,  $c, k$ remain continuous.
\end{itemize}

Thanks to these wave properties, 
the solution of the Riemann problem for~\eqref{1.1}-\eqref{1.3} 
is rather simple. 
Given the left and right states $(s_l,c_l,k_l)$ and $(s_r,c_r,k_r)$,
we have the following  global Riemann solver:
\begin{itemize}
\item Let $(s_m,c_l, k_r)$  denote the right state of the $k$-wave. The value  $s_m$
is uniquely determined by the condition 
$$ f(s_m,c_l,k_r) = f(s_l,c_l,k_l).$$

\item For the remaining waves, we have $k\equiv k_r$ throughout. 
We then solve the Riemann problem for the $2\times 2$ system 
\eqref{1.1}-\eqref{1.2} with Riemann data 
$(s_m,c_l)$ and $(s_r,c_r)$ as the left and right states.
We use the Riemann solver in section~\ref{sec2.1}.
The solution consists of waves with non-negative speed.
\end{itemize}

%%%%%%%%%%%%%%%%%%
\begin{remark}\label{rm1}
It would be of interest to prove existence and well-posedness result for the Cauchy problem
directly for the Eulerian system~\eqref{1.1}-\eqref{1.3}. 
The key estimate is the bound on the total wave strength, suitably defined.
We expect that
the resonance between the $s$~and $c$~families can be controlled by 
Temple-style functionals for wave strength.
Then,  the total wave strength is non-increasing at interactions between
$s$ and $c$ waves. 
However, there are additional difficulties caused by
the interactions between $s$ and $k$~waves.
In strictly hyperbolic cases, this can be controlled by
adding a suitable interaction potential functional.
However, the interaction potential functional here must take into account 
the difficulties caused by the vacuum state where $s=0$.
\end{remark}

%%%%%%%%%%%%%%%%%%%%%%%%%%%
\section{Polymer flooding with adsorption in rough media}\label{sec3}
\setcounter{equation}{0}
%%%%%%%%%%%%%%%%%%%%%%%%%%%

In this section 
we consider the polymer flooding with
the adsorption effect, but neglect the gravity force:
\begin{align}
 s_t + f(s,c,k)_x =0, \label{2.1} \\
 (m(c)+cs)_t + (cf(s,c,k))_x =0,\label{2.2} \\
 k_t =0.\label{2.3} 
\end{align}
The new term $m(c)$ denotes the adsorption of polymers into the porous media,
which we assume to be a function of $c$.
Typically one assumes that 
$$ m'(c)>0, \qquad m''(c)<0, \qquad \forall c. $$

The same coordinate change as \eqref{eq:c} leads to the following Lagrangian system:
\begin{align}
 \left(\frac{s}{f}\right)_\phi -\left(\frac{1}{f}\right)_\psi =0 ,\label{2.4}\\
 m(c)_\phi + c_\psi =0,\label{2.5}\\
k_\phi=0.\label{2.6}
\end{align} 
The equation~\eqref{2.6} for $k$ is unchanged. 
The equation~\eqref{2.5} for $c$ is stilled decoupled, but now
$c$  solves a scalar conservation law, and
the $c$~family is  genuinely nonlinear. 
Nevertheless, $c$ can be solved independently. 
Given initial data for $k,c$, their values at any point $(\phi,\psi)$
can be computed first. Then,  $s$ solves the scalar conservation law~\eqref{eqL}
with discontinuous coefficient.
The discontinuities include the jumps in $k$ and shocks in the solution of  $c$, 
which has a more complex structure, 
see Remark~\ref{rm2}.

%%%%%%%%%%%%%%%%
\subsection{The reduced $2\times2$ model}\label{sec3.1}
%%%%%%%%%%%%%%%%

When $k$ is constant, one has the reduced $2\times2$ system
\begin{equation}\label{eq:red}
\begin{cases} s_t + f(s,c)_x =0, \\ 
 (m(c)+cs)_t + (cf(s,c))_x =0.\end{cases}
\end{equation}
Given Riemann data $(s_l,c_l),(s_r,c_r)$,
the Riemann problem is studied in the literature, even for multi-component polymers,  
see~\cite{JW1, MR1000728, Dahl, JTW}.
In particular, for any $(s_l,s_r)$, 
if $c_l>c_r$, the solution contains a $c$~shock; 
and if $c_l<c_r$ then we have a $c$~rarefaction. 

Consider the case $c_l>c_r$ where we have a $c$~shock.
We define
\begin{equation}\label{agg}
 a\,\dot=\,\frac{m(c_l)-m(c_r)}{c_l-c_r},\qquad g_l(s) ~\dot=~\frac{f(s,c_l)}{s+a}, \qquad
 g_r(s) ~\dot=~\frac{f(s,c_r)}{s+a}.
 \end{equation}
Define the functions $G^\flat,G^\sharp$ as in \eqref{Gsharp}-\eqref{Gflat},
and let $s_\pm$ be the {\em minimum jump path} that connects $G^\sharp$ and $G^\flat$.
Then, $s_\pm$ will be the trace along $c$ shock which travels with speed
$$\sigma = g_l(s_-)=g_r(s_+).$$
Once the $c$~shock is located, we patch up the solutions of a regular conservation 
law on the left and on the right of the $c$-shock, as described in section~\ref{sec2.1}.

\begin{figure}[htbp]
\begin{center}
\includegraphics[ width=7cm,clip, trim=1mm 7mm 1mm 11mm]{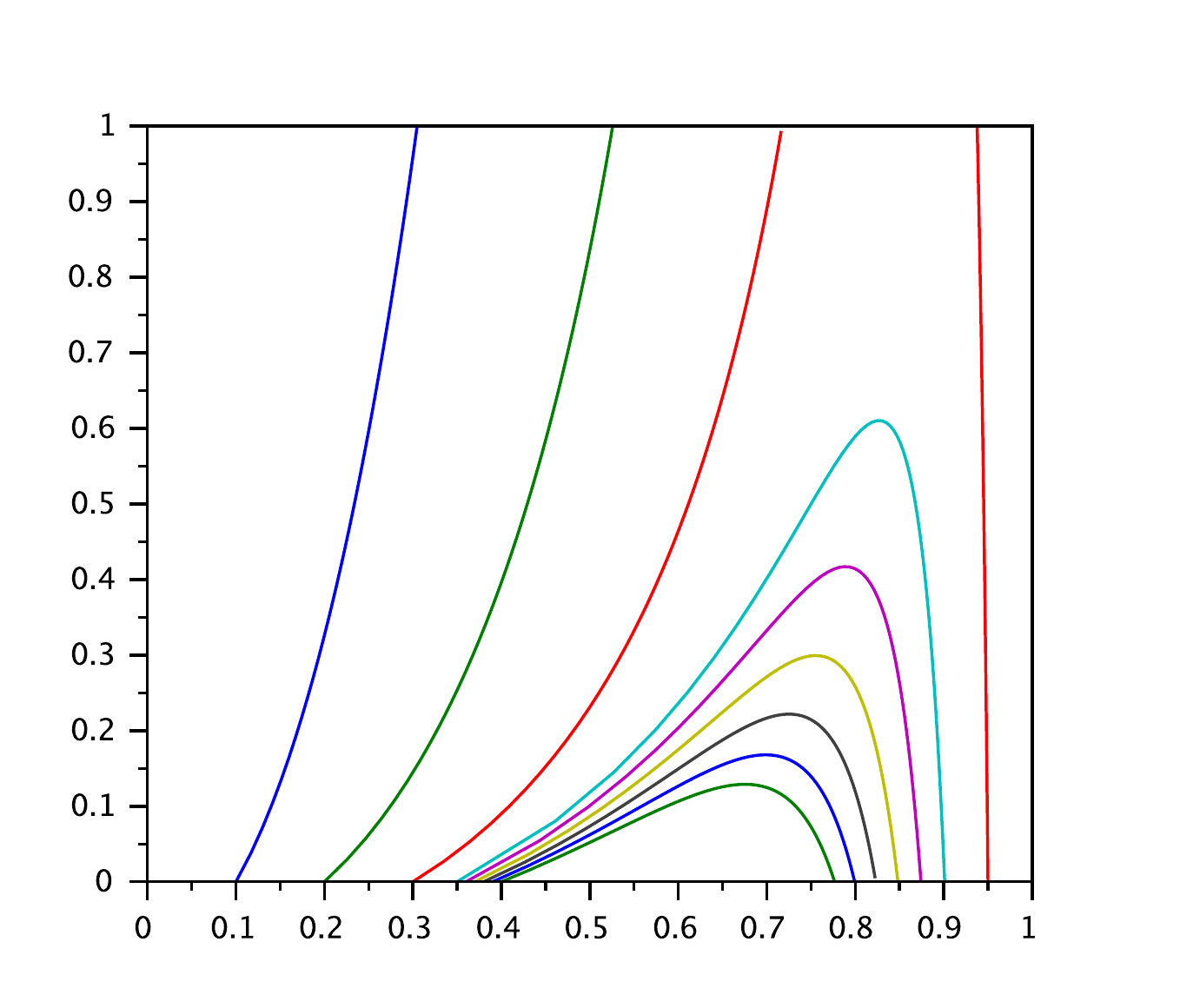}
\caption{Integral curve for $c$ family in the $(s,c)$-plane. }
\label{fig:r2}
\end{center}
\end{figure}

When $c_l<c_r$ and we have a $c$ rarefaction wave, 
the path of the rarefaction wave goes along the integral curves of the $c$~eigenvectors.
See Figure~\ref{fig:r2} for a typical graph of these integral curves. 
The resonance point occurs at where the integral curves have horizontal tangent.
Thus, the $c$ rarefaction path must lie either on the left or on the right of the 
resonance point. 
Unique path can be chosen to allow feasible solution with increasing wave
speeds from left to right.  See for example~\cite{JW1}. 
We omit further details.

%%%%%%%%%%%%%%%%
\subsection{Wave properties and a global Riemann solver} 
%%%%%%%%%%%%%%%

We go back to the  full $3\times3$ system~\eqref{2.1}-\eqref{2.3}. 
The Jacobian matrix for the flux function is again triangular
$$
J=\begin{pmatrix}
f_s & f_c & f_k \\
0 & f/(s+m'(c)) & 0 \\
0&0&0
\end{pmatrix}
$$
with three eigenvalues
$$
\lambda_s = f_s(s,c,k), \qquad 
\lambda_c =\frac{f(s,c,k)}{s+m'(c)} , \qquad
\lambda_k=0,
$$
and three  corresponding right-eigenvectors
$$
r_s=\begin{pmatrix} 1 \\ 0 \\ 0\end{pmatrix},
\qquad
r_c=\begin{pmatrix} -f_c(s,c,k) \\ \lambda_s-\lambda_c \\ 0 \end{pmatrix},
\qquad
r_k=\begin{pmatrix}  -f_k(s,c,k) \\ 0 \\ f_s(s,c,k)\end{pmatrix}.
$$
Direct computations give the following directional derivatives:
\begin{align}
\nabla \lambda_s \cdot r_s &~=~ f_{ss}(s,c,k), \nonumber \\
\nabla \lambda_c \cdot r_c &~=~ \frac{-f(s,c,k) m''(c) (\lambda_s-\lambda_c)}{(s+m'(c))^2}, 
\label{eig}
\\
\nabla \lambda_k \cdot r_k &~=~0.\nonumber
\end{align}

We observe the following properties.
\begin{itemize}
\item
The $k$ family is linearly degenerate and travels with speed 0, which is the slowest family.
Crossing a $k$ contact, both $f$ and $c$ remain continuous.
\item
The $c$ family is genuinely nonlinear, as indicated by \eqref{2.5} in the Lagrangian system.
We can have either a single $c$ shock or a single
$c$ rarefaction fan in the solution.  

This fact is not clear from the directional derivatives~\eqref{eig} in the Eulerian system.
When $\lambda_s\not=\lambda_c$,  we can rewrite the $c$~eigenvector as
$$
\tilde r_c=\begin{pmatrix} \displaystyle \frac{f_c(s,c,k)}{ \lambda_c-\lambda_s} \\[3mm]  1 \\[2mm]  0 \end{pmatrix}.
$$
The integral curves for $\tilde r_c$ is now parametrized by $c$.  
Straight computations give
$$
\nabla \lambda_c \cdot \tilde r_c = -\frac{f(s,c,k) m''(c)}{(s+m'(c))^2} >0.
$$
Therefore, we have a $c$ rarefaction when $c_l<c_r$, and the rarefaction curve
can never cross the resonant point where $\lambda_s=\lambda_c$.
When $c_l>c_r$, we have a $c$ shock.
\item
The $s$ family is not genuinely nonlinear, but it's a Temple family where shock curves and 
rarefaction curves coincide. Crossing a $s$ wave, both $c$ and $k$ remain continuous. 
\end{itemize}

Given any Riemann data $(s_l,c_l,k_l), (s_r,c_r,k_r)$, 
we now have the following global Riemann solver, similar to section~\ref{sec2.3}:
\begin{itemize}
\item
The right state of the $k$ wave is $(s_m,c_l,k_r)$ where $s_m$ is uniquely determined 
by 
$$f(s_l,c_l,k_l)=f(s_m,c_l,k_r).$$
\item
For the remaining waves, we have $k\equiv k_r$ which is constant, 
Then, we solve the Riemann problem for the reduced model~\eqref{eq:red},
with left and right states $(s_m,c_l)$ and $(s_r,c_r)$ respectively, 
following the Riemann solver in section~\ref{sec3.1}.
\end{itemize}

\begin{remark}\label{rm2}
We remark that wave interaction estimates for this system remain very complicated,
and the control of the total wave strength is not available in the literature. 
In a recent work~\cite{GS-preprint}, a scalar conservation law
with general discontinuous  flux is studied
$$ u_t + f(\alpha(t,x),u)_x =0, \qquad u(0,x)=\bar u(x).$$
Here $\alpha(t,x)$ is discontinuous w.r.t.~both variables $t,x$.
Recall that a function of a single variable $\alpha: I\!\!R \mapsto I\!\!R$
is {\em regulated} if it admits left and right limits at every point. 
Such a concept can be extend to functions of two variables. 
Assuming that $\alpha(t,x)$ is regulated, in~\cite{GS-preprint}  we prove
that the vanishing viscosity solutions of 
$$ u_t + f(\alpha(t,x),u)_x =\ve u_{xx}, \qquad u(0,x)=\bar u(x)$$
converge to a unique limit solution. 

One can show that solutions of scalar conservation laws with
convex flux function are  regulated functions.
An extension of the result in~\cite{GS-preprint} could prove
a similar result, at least for the Lagrangian system~\eqref{2.4}-\eqref{2.6}. 
Details may come in future works.
\end{remark}

%%%%%%%%%%%%%%%%%%%%%%%%%%%
\section{A second order traffic flow model with rough road condition}\label{sec4}
\setcounter{equation}{0}
%%%%%%%%%%%%%%%%%%%%%%%%%%%

As a model of intermediate level of complexity, 
we consider a $3\times3$ system for traffic flow
\begin{align}
\rho_t + (\rho v)_x &=0, \label{eq:1}\\
[\rho(v+k \rho^\gamma)]_t + [\rho v(v+k \rho^\gamma)]_x &=0,\label{eq:2}\\
k_t &=0.\label{eq:3}
\end{align}
Here $\rho \ge 0$ denotes the car density, $v \ge 0$ is car velocity, 
and $k(x) >0$ denotes the road condition. 
Furthermore, $\gamma\in(1,2)$ is a constant.
We consider rough road condition, where $k(x)$ is discontinuous.

When $k$ is constant, the reduced system \eqref{eq:1}-\eqref{eq:2} 
was proposed in~\cite{AR}. 
Equation~\eqref{eq:1} denotes the conservation of mass. 
In~\eqref{eq:2},  the quantity $ k  \rho^\gamma$ denotes some kind of ``pressure''. 
The  physical modeling  leads to  a non-conservative 
formulation
\begin{equation}\label{eq:VP}
(v+k  \rho^\gamma)_t + v \cdot (v+k  \rho^\gamma)_x=0.
\end{equation}
With some algebraic manipulation and utilizing \eqref{eq:1}, 
one can rewrite \eqref{eq:VP} in the conservative form of~\eqref{eq:2}.  
Although equation~\eqref{eq:2} resembles the conservation of momentum, 
there is no physical meaning for the conserved quantity $\rho (v+k  \rho^\gamma)$. 

For notational convenience, we denote that 
\begin{equation}\label{eq:defw}
w~\dot=~v+ k  \rho^\gamma.
\end{equation}
Note that if $w$=constant and $\rho>0$, 
then \eqref{eq:2} reduces to \eqref{eq:1}.

%%%%%%%%%%%%%%%%%%%%%%%%%%%
\subsection{A Lagrangian system and the decoupling feature}
%%%%%%%%%%%%%%%%%%%%%%%%%%%
Consider  a Lagrangian coordinate $(\phi,\psi)$ defined as 
$$ 
\phi_x = -\rho, \quad \phi_t = \rho v, \qquad \phi(0,0)=0,\qquad
\psi=x.
$$
When $\rho v >0$, the coordinate change is well-defined.
Direct computation leads to the following Lagrangian system
\begin{align}
\left( \frac{1}{\rho v}\right)_\psi - \left( \frac{1}{v}\right)_\phi  &= 0, \label{b1.4}\\
w_\psi &=0,\label{b1.5}\\
 k _\phi &=0.\label{b1.6}
\end{align}

We observe the decoupling feature for $k$ and $w$  in this Lagrange coordinate.
Herr $ k  $ is contant in $\phi$, 
and  $w$  is constant in $\psi$.
These features are very similar to those of the system~\eqref{1.6}-\eqref{1.8}.
Given the initial data at $t=0$, the values of $(w,k)$ for any 
coordinate point $(\phi,\psi)$ are determined 
trivially,  see Figure~\ref{fig:1}.  
Once $( k (\phi,\psi) , w(\phi,\psi))$  are given, 
we can express $v$ in terms of $\rho$, 
i.e., 
\begin{equation}\label{eq:v}
 v= v(\rho; w, k )= 
w- k  \rho^\gamma.
\end{equation}
Then it remains to solve $\rho$ using the scalar conservation law 
\eqref{b1.4} with variable coefficient, 
\begin{equation}\label{LL}
\left(\frac{1}{\rho \cdot (w- k  \rho^\gamma)} \right) _\psi -
\left(\frac{1}{ w- k  \rho^\gamma} \right)_\phi =0,
\end{equation}
where $(w,k)$ are given functions, possibly discontinuous, and $\rho$ is
the unknown.
The discontinuities in the flux function occur along 
horizontal and vertical lines in the $(\phi,\psi)$ coordinate. 

It would be interesting to explore possible ways of extending 
the result in~\cite{CocliteRisebro} 
for this case, taking extra care of the vacuum state.

%%%%%%%%%%%%%%%%%%%%%%%%%%%
\subsection{Some basic analysis}
%%%%%%%%%%%%%%%%%%%%%%%%%%%

For the Eulerian system~\eqref{eq:1}-\eqref{eq:3},
treating $(\rho, v,k)$ as the unknown vector, 
the Jacobian matrix for the flux function  is triangular
$$
J=\begin{pmatrix}
v & \rho & 0 \\
0 & v-\gamma  k  \rho^{\gamma} & v \rho^\gamma\\
0&0&0
\end{pmatrix}.
$$
Denoting the three families as the $\rho$, $v$, and $k$ families, 
we have three eigenvalues
$$ \lambda_\rho = v, \qquad \lambda_v = v-\gamma  k  \rho^{\gamma} ,
\qquad \lambda_k =0,
$$
with three corresponding right-eigenvectors
$$
r_\rho =\begin{pmatrix} 1\\ 0\\ 0\end{pmatrix},
\qquad
 r_v=\begin{pmatrix} -1 \\ \gamma  k  \rho^{\gamma-1} \\ 0\end{pmatrix},
\qquad
 r_ k =\begin{pmatrix} \rho^{\gamma+1}\\ -v \rho^\gamma \\ 
 v-\gamma  k  \rho^{\gamma} \end{pmatrix}.
$$

Straight computations give 
\begin{eqnarray*}
\nabla \lambda_\rho \cdot r_\rho &=& 0, \\
\nabla \lambda_v \cdot r_v
&=& (\gamma^2+\gamma)  k  \rho^{\gamma-1} >0,\\
\nabla \lambda_k \cdot r_k &=& 0.
\end{eqnarray*}

Thus, the $\rho$ and $k$ families are linearly degenerate,
where all discontinuities are contacts. 
The $v$ family is genuinely nonlinear, where we have either a  $v$ shock or a
$v$ rarefaction in the solution of a Riemann problem. 
The $v$ rarefaction curves are integral curves of $r_v$. 

We consider the $\rho$ and $v$ jumps. 
Observe that crossing both $\rho$ and  $v$ waves, the value $ k $ remains constant.  
Fix a $ k $ value, we consider a jump initiated from $(\rho_0,v_0)$. 
Writing $w_0= v_0+ k  \rho_0^\gamma$, and letting  $\sigma$ be the jump speed,
the RH jump conditions require
\begin{align}
\sigma (\rho - \rho_0) &= \rho v - \rho_0 v_0, \label{eq:a}\\
\sigma (\rho w-\rho_0 w_0) &= \rho v w - \rho_0 v_0 w_0.\label{eq:a2}
\end{align}

We first consider the case with vacuum. 
If $\rho_0=0$, then for any values of $ v_0$  we have 
$$
\mbox{either} ~ 
\Big\{\rho=0, ~ \sigma~\mbox{arbitrary}\Big\}
\qquad \mbox{or}\qquad
\Big\{\sigma=v, ~ \rho~\mbox{arbitrary} \Big\}.
$$
On the other hand, if $\rho=0$,  then for any values of $v$, we have 
$$
\mbox{either} ~ 
\Big\{\rho_0=0, ~ \sigma~\mbox{arbitrary}\Big\}
\qquad \mbox{or}\qquad
\Big\{\sigma=v, ~ \rho_0~\mbox{arbitrary} \Big\}.
$$
We remark that the vacuum state is special,
where the $v$ and $\rho$ families
have the same eigenvalue and eigenvector, 
so  the system is both parabolic degenerate and 
linearly degenerate. 

\medskip

For the rest, we assume $\rho,\rho_0>0$. 
If $v=v_0$, \eqref{eq:a} gives $\sigma=v_0=v$, and \eqref{eq:a2}  trivially holds. 
This gives a $\rho$ contact discontinuity.  
Note that $\sigma=v_0$ or $\sigma=v$ leads to the same wave. 
We conclude that, crossing a $\rho$~wave,   both $ k , v$ remain continuous.

Now we consider the $v$~shocks,  assuming 
$v\not=v_0\not=\sigma$.
Multiplying \eqref{eq:a} by $w$ and subtracting from \eqref{eq:a2}, 
we get 
$$ \rho_0 (\sigma - v_0) (w-w_0)= 0.$$
Or symmetrically, 
multiplying \eqref{eq:a} by $w_0$ and subtracting from \eqref{eq:a2},
we get
$$
 \rho(\sigma - v) ( w-w_0)= 0.
$$
Since $\rho,\rho_0>0$, this implies 
$$ w=w_0.$$
The $v$ shock travels with speed:
$$
\sigma_v= \frac{\rho v -\rho_0 v_0}{\rho -\rho_0}.
$$

We further observe that, along a $v$~integral curve, $w$  remains constant. 
Indeed, we have
$$
\nabla w \cdot r_v = 
\begin{pmatrix} \gamma  k  \rho^{\gamma-1} \\ 1 \\  \rho^\gamma\end{pmatrix}
\cdot \begin{pmatrix} -1 \\ \gamma  k  \rho^{\gamma-1} \\ 0\end{pmatrix}
= 0.
$$
We conclude that the $v$ rarefaction curves coincide with $v$ shock curves.
Thus, the $3\times3$ system~\eqref{eq:1}-\eqref{eq:3} is a Temple class, 
where the system is of mixed type.

We remark that, 
$$ \lambda_\rho \ge \lambda_v, \qquad \lambda_\rho \ge \lambda_k,$$
but $\lambda_v-\lambda_k$ may change sign.
Thus the possible nonlinear resonance  only occurs between 
the linearly degenerate $k$~family and the genuinely nonlinear $v$~family.
This fact  should make it possible to control the  resonance.

\medskip 

We  summarize the wave behaviors:  
\begin{itemize}
\item 
Crossing a $ k $-contact, both $\rho v$ and $w$ remain continuous.
\item
Crossing a $v$-front, $w$ and $ k $ remain continuous. 
\item
Crossing a $\rho$-contact, $ k $ and $v$ remain continuous.  
\end{itemize}
Vacuum state is considered as a mixing of $v$ and $\rho$ fronts,
as it will be clear later. 
Note also that $(w,v)$ serve as the natural Riemann invariants for the 
$(v,k)$ families, respectively.

%%%%%%%%%%%%%%%%%%%%%%%%%%%
\subsection{The reduced model and its global Riemann solver}\label{sec4.2}
%%%%%%%%%%%%%%%%%%%%%%%%%%%

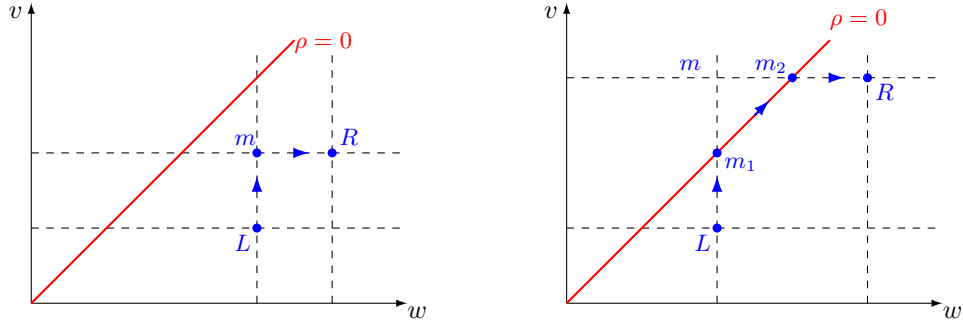
\begin{figure}[htbp]
\begin{center}
\setlength{\unitlength}{1mm}
\begin{picture}(70,40)(0,0)  % -- Left ---
\put(0,0){\vector(1,0){50}}\put(50,-2){$w$}
\put(0,0){\vector(0,1){40}}\put(-3,38){$v$}
\multiput(0,10)(2,0){25}{\line(1,0){1}}
\multiput(0,20)(2,0){25}{\line(1,0){1}}
\multiput(30,0)(0,2){17}{\line(0,1){1}}
\multiput(40,0)(0,2){17}{\line(0,1){1}}
\thicklines
\color{red}
\put(0,0){\line(1,1){35}} \put(35, 34){\small $\rho=0$}
\color{blue}
\put(30,10){\circle*{1.2}}\put(30,20){\circle*{1.2}}\put(40,20){\circle*{1.2}}
\put(27,7){\small $L$} \put(27,21){\small $m$} \put(41,21){\small $R$}
\put(30,15){\vector(0,1){2}}\put(35,20){\vector(1,0){2}}
\end{picture}
% -------------------------------
\begin{picture}(50,40)(0,0)  % -- Right ---
\put(0,0){\vector(1,0){50}}\put(50,-2){$w$}
\put(0,0){\vector(0,1){40}}\put(-3,38){$v$}
\multiput(0,10)(2,0){25}{\line(1,0){1}}
\multiput(0,30)(2,0){25}{\line(1,0){1}}
\multiput(20,0)(0,2){17}{\line(0,1){1}}
\multiput(40,0)(0,2){17}{\line(0,1){1}}
\thicklines
\color{red}
\put(0,0){\line(1,1){35}} \put(35, 37){\small $\rho=0$}
\color{blue}
\put(20,10){\circle*{1.2}}\put(20,20){\circle*{1.2}}
\put(30,30){\circle*{1.2}}\put(40,30){\circle*{1.2}}
\put(17,7){\small $L$} \put(21,18){\small $m_1$} 
\put(25,31){\small $m_2$} \put(41,27){\small $R$}
\put(15,31){\small $m$}
\put(20,15){\vector(0,1){2}}\put(25,25){\vector(1,1){2}}\put(35,30){\vector(1,0){2}}
\end{picture}
\caption{Riemann solver for the $2\times2$ model. Left: Without vacuum, the solution consists of a $v$-rarefaction and a $\rho$-contact. 
 Right: With an intermediate vacuum wave.}
\label{fig:RP2}
\end{center}
\end{figure}

When $k$  is constant, say $k=1$,  \eqref{eq:1}-\eqref{eq:3} reduces to 
a $2\times 2$ system
\begin{equation}\label{ar}
\begin{cases} 
\rho_t + (\rho v)_x =0, \\ 
\left[\rho(v+ \rho^\gamma)\right]_t + \left[\rho v(v+ \rho^\gamma)\right]_x =0.
\end{cases}
\end{equation}
The Riemann problem for~\eqref{ar} was studied in much detail in \cite{AR}.
However, utilizing the Riemann invariants $(w,v)$, 
the Riemann solver can be presented in a very compact manner,
as illustrated in Figure~\ref{fig:RP2}.
Note that the vacuum state $\rho=0$ is the straight line $w=v$  in the $(w,v)$-plane,
indicated by the red line. 
The physically feasible region  $\rho\ge 0$ lies on the right side of the vacuum line.
Let $L=(w_l,v_l)$ and $R=(w_r,v_r)$ be the Riemann data, both lie on the right of the vacuum line.
Consider the state $m=(w_l,v_r)$. We have two cases:
\begin{itemize}
\item 
If $m$ lies on the right side of the vacuum line, then the 
solution of the Riemann problem consists of two waves, 
with a $v$-wave connecting $L$ to $m$, followed by a $\rho$-wave connecting
$m$ to $R$. 
The $v$ wave is a shock if  $v_l>v_r$, and a rarefaction if $v_l<v_r$.
See the left plot in Figure~\ref{fig:RP2}.
\item 
If $m$ lies on the left of the vacuum line, then vacuum occurs in the solution,
see the right plot in Figure~\ref{fig:RP2}. 
This could only happen when $v_l<v_r$. 
We have two intermediate states
$m_1,m_2$ where $\rho=0$. 
The solution of the Riemann problem consists of 3 waves:
a $v$ rarefaction connecting $L$ to $m_1$, 
followed by a vacuum wave connecting $m_1$ to $m_2$, 
and finally a $\rho$ contact from $m_2$ to $R$.
\end{itemize}

The vacuum wave is rather ``fake'', since $\rho\equiv 0$ is always a solution for any values of $v$.
To ``assign'' the $v$ values along a vacuum wave in the solution of a Riemann problem, 
we set $\rho=0$ in \eqref{eq:2} and obtain the 
Burgers' equation 
$$ v_t + (v^2/2)_x=0.$$ 
Since $v$ increases from $m_1$ to $m_2$, 
the solution for $v$ is a rarefaction wave.

\textbf{Continuous dependence.} 
Viewed in the $(w,v)$ plane, it is clear that
the path for the solution of a Riemann problem depends 
continuously on the data $L$ and $R$.

\begin{figure}[htbp]
\begin{center}
\setlength{\unitlength}{0.8mm}
\begin{picture}(70,40)(0,0)  % -- Left -- 
\multiput(0,10)(2,0){25}{\line(1,0){1}}
\multiput(0,30)(2,0){25}{\line(1,0){1}}
\multiput(40,0)(0,2){17}{\line(0,1){1}}
\thicklines
\color{red}
\put(0,0){\line(1,1){35}} \put(35, 37){\small $\rho=0$}
\color{blue}
\put(10,10){\circle*{1.2}}
\put(40,10){\circle*{1.2}}
\put(30,30){\circle*{1.2}}
\put(37,6){\small $L$}  \put(24,31){\small $R$} \put(10,6){\small $M$}
\color{cyan}
\put(40,30){\circle*{1.2}}\put(41,31){\small $m$}
\end{picture}
\begin{picture}(50,40)(0,-3)  % -- Right -- 
\multiput(0,20)(2,0){25}{\line(1,0){1}}
\multiput(0,2)(2,0){25}{\line(1,0){1}}
\multiput(20,0)(0,2){17}{\line(0,1){1}}
\multiput(30,0)(0,2){17}{\line(0,1){1}}
\thicklines
\color{red}
\put(0,0){\line(1,1){35}} \put(35, 32){\small $\rho=0$}
\color{blue}
\put(20,20){\circle*{1.2}}
\put(30,30){\circle*{1.2}}\put(30,2){\circle*{1.2}}
\put(16,21){\small $L$}  \put(24,31){\small $M$} \put(31,-2){\small $R$}
\color{cyan}
\put(20,2){\circle*{1.2}}\put(15,-1){\small $m$}
\end{picture}
\caption{Left: interaction of a vacuum wave with a $\rho$ wave.
Right: interaction of a vacuum wave with a $v$ wave.
Here $M$ is the middle state of the incoming waves, and $m$ 
is the middle state of the outgoing waves.}
\label{fig:RP2ex}
\end{center}
\end{figure}
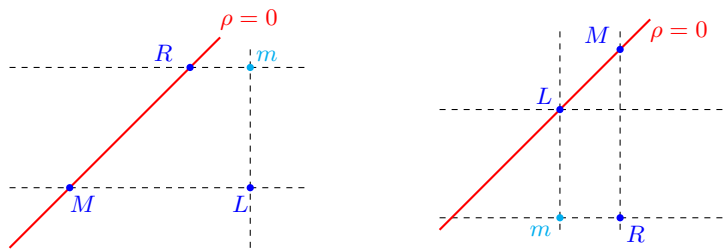

\textbf{Interaction estimates.} 
One can define the wave strength as the a Manhattan distance in the $(w,v)$-plane,
i.e., 
any wave connecting $(w_l,v_l)$ and $(w_r,v_r)$  has the strength
$$ \abs{w_l-w_r}+\abs{v_l-v_r}.$$
We claim that the total wave strength remains non-increasing at any interaction.
Indeed,  we have the following observations:
\begin{itemize}
\item 
Two $\rho$~waves can not interact with each other since the family is linearly 
degenerate.
\item
When a $\rho$~wave interacts with a $v$~wave, the total wave strength is unchanged. 
\item
For interactions between two $v$~waves, the total wave strength is non-increasing since 
the family is genuinely nonlinear.
\item
When a  vacuum wave interacts with either a $v$~wave or a $\rho$~wave, 
cancellation happens and the total wave strength  is decreasing. 
See Figure~\ref{fig:RP2ex}. 
\end{itemize}

\textbf{Front tracking.} 
In a front tracking approximation, we approximate $v$~rarefaction waves
with upward jumps of size $\ve$. 
It's simple to show that the algorithm is well-posed.
Total number of fronts is uniformly bounded. 
Total wave strength, measured with the Manhattan distance, 
 is also uniformly bounded. 
Existence of entropy weak solution for the Cauchy problem follows
from standard theory.

 \medskip
 
\textbf{Critics for the model.}
This second order traffic flow model admits some unreasonable solutions.
For example, if $v(0,x)\equiv 0$, then $\rho_t=0$ and 
we have the solution $\rho(t,x)=\rho(0,x)$ for all $t>0$. 
This means, if cars are initially stationary on a road, they will remain
stationary for all time. This unreasonable behavior is 
caused by the conservation of the  ``momentum'' $\rho(v+ k  \rho^\gamma)$,
a concept borrowed from gas dynamics. 
However, 
moving cars behave differently from gas particles,
and the momentum should not be conserved.
High order models for traffic flow are better formulated with a relaxation parameter, 
where one considers the reaction/acceleration time for each driver.

%%%%%%%%%%%%%%%%%%%%%%%%
\subsection{Riemann solver for the $3\times3$ system}
%%%%%%%%%%%%%%%%%%%%%%%%

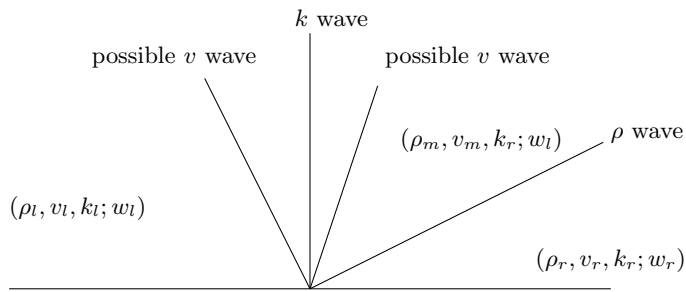
\begin{figure}[htbp]
\begin{center}
\setlength{\unitlength}{1mm}
\begin{picture}(80,40)(-40,0)  
\put(-40,0){\line(1,0){80}}
\put(0,0){\line(2,1){39}}\put(40,20){\small $\rho$ wave}
\put(0,0){\line(0,1){34}}\put(-2,35){\small $k$ wave}
\put(0,0){\line(1,3){9}}\put(10,30){\small possible $v$ wave}
\put(0,0){\line(-1,2){14}}\put(-29,30){\small possible $v$ wave}
\put(-40,10){\small $(\rho_l,v_l,k_l;w_l)$}
\put(30,3){\small $(\rho_r,v_r,k_r;w_r)$}
\put(12,19){\small $(\rho_m,v_m,k_r;w_l)$}
\end{picture}
\caption{Wave structures  in the Riemann solution for the $3\times3$ traffic flow model.}
\label{fig:ww}
\end{center}
\end{figure}

We now describe a global Riemann solver for~\eqref{eq:1}-\eqref{eq:3}. 
Let $(\rho_l,v_l,k_l), (\rho_r,v_r,k_r)$ denote the Riemann data,
and $w_l,w_r$ be the correspond $w$ values. 
Since the $\rho$~wave is the fastest one,
it will have $(\rho_r,v_r,k_r)$ as its right state. 
Denote the left state of the $\rho$~wave by $(\rho_m,v_m,k_r)$. 
Since $w$ is constant crossing both $k$  and $v$ waves, 
we have $w_l$ on the left of the $\rho$~wave.  
See Figure~\ref{fig:ww}. 

The global Riemann solver consists of two steps. 

\begin{figure}[htbp]
\begin{center}
\setlength{\unitlength}{0.8mm}
\begin{picture}(70,42)(0,-3)  % -- Left ---
\put(0,0){\vector(1,0){50}}\put(50,-2){$w$}
\put(0,0){\vector(0,1){40}}\put(-3,38){$v$}
\multiput(0,20)(2,0){25}{\line(1,0){1}}
\multiput(25,0)(0,2){12}{\line(0,1){1}}
\multiput(35,0)(0,2){12}{\line(0,1){1}}
\thicklines
\color{red}
\put(0,0){\line(1,1){35}} \put(36, 34){\small $\rho=0$}
\color{blue}
\put(25,20){\circle*{1.2}}\put(35,20){\circle*{1.2}}
\put(26,21){\small m} \put(36,21){\small R}
\put(24,-3){\small $w_l$} \put(34,-3){\small $w_r$}
\put(-5,19){\small $v_r$}
\end{picture}
% -------------------------------
\begin{picture}(55,42)(-5,-3)  % -- Right ---
\put(0,0){\vector(1,0){50}}\put(50,-2){$w$}
\put(0,0){\vector(0,1){40}}\put(-3,38){$v$}
\multiput(0,25)(2,0){25}{\line(1,0){1}}
\multiput(15,0)(0,2){10}{\line(0,1){1}}
\multiput(35,0)(0,2){14}{\line(0,1){1}}
\multiput(0,15)(2,0){10}{\line(1,0){1}}
\thicklines
\color{red}
\put(0,0){\line(1,1){35}} \put(37, 34){\small $\rho=0$}
\color{blue}
\put(15,15){\circle*{1.2}}\put(25,25){\circle*{1.2}}\put(35,25){\circle*{1.2}}
\put(16,12){\small m$_1$} \put(25,22){\small m$_2$} \put(36,26){\small R}
\put(14,-3){\small $w_l$}\put(34,-3){\small $w_r$}
\put(-5,25){\small $v_r$}\put(-5,14){\small $v_m$}
\end{picture}
\caption{Riemann solver for the $3\times3$ model: The algorithm that determines
the $\rho$-front. Left: No vacuum.  Right: With vacuum and a  $v$-rarefaction attached
on the left of the $\rho$-front.}
\label{fig:RP3x}
\end{center}
\end{figure}
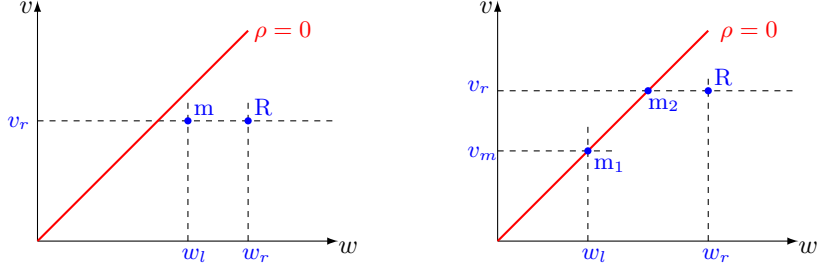

\textbf{Step 1.} 
We determine the value $(\rho_m,v_m)$.
There are two cases, with and without the vacuum state.
\begin{itemize}
\item
 If  $ w_l \ge v_r$ (see left plot in Figure~\ref{fig:RP3x}), 
 we can compute the unique value of $\rho_m$ using
 $$ v_r +  k _r \rho_m^\gamma = w_l. $$
 This gives 
 $$
 \rho_m = \left(\frac{w_l-v_r}{ k _r}\right)^{1/\gamma} , 
 \qquad
 v_m = v_r.
 $$
 \item
 Otherwise if $w_l<v_r$, we have a vacuum wave in the solution
 (see the right plot in Figure \ref{fig:RP3x}). 
 From $m_2$ to $R$ we have a $\rho$ wave.
 On its left there is a vacuum wave that connects $m_1$ to $m_2$. 
 Denote the left state of the vacuum wave as $(\rho_m,v_m,k_r)$, we set 
 $$\rho_m=0, \qquad v_m=w_l.$$ 
 \end{itemize}

\textbf{Step 2.}  As the second step, 
one solves a Riemann problem for the two states
$$ (\rho_l, v_l,  k _l; ~w_l) , \qquad (\rho_m, v_m,  k _r;  ~w_l)$$
with only $k$ wave and $v$ waves. 
Since $w\equiv w_l$ throughout the solution, \eqref{eq:2} reduces to \eqref{eq:1}.
Furthermore, we also have
$$v=v(\rho; k )=w_l- k  \rho^\gamma.$$
It remains to solve a scalar conservation law with discontinuous flux: 
\begin{align}
 \rho_t + \left(\rho(w_l-  k (x) \rho^\gamma )\right)_x =0,
 \end{align}
where $ k (x)$ is the jump function connecting $ k _l, k _r$ at $x=0$.
Calling the flux functions
$$
f_l(\rho)=g_l(\rho) = \rho(w_l-  k_l \rho^\gamma ), 
\qquad
f_r(\rho)=g_r(\rho) = \rho(w_l-  k_r \rho^\gamma ),
$$
and defining $G^\sharp,G^\flat$ accordingly as in \eqref{Gsharp}-\eqref{Gflat}, 
replacing $s$ with $\rho$.
Then, the $k$~wave is located at the {\em minimum jump path}
that connects $G^\sharp$ and $G^\flat$. 
The remaining waves in this Riemann solver are 
determined by patching up solutions, as in~\eqref{Rl}-\eqref{Rr}.

\bigskip

In conclusion, we  have constructed 
a global Riemann solver that generates a unique self similar solution 
for any given left and right state.  In the solution, all the quantities 
$(\rho, v,  k , w)$ are  non-negative.

%%%%%%%%%%%%%%%%%%%%%%%%%%%
\section{Polymer flooding with gravity and rough media}\label{sec5}
\setcounter{equation}{0}
%%%%%%%%%%%%%%%%%%%%%%%%%%%

\begin{figure}[htbp]
\begin{center}
\includegraphics[width=7cm,clip,trim=0mm 4mm 5mm 10mm]{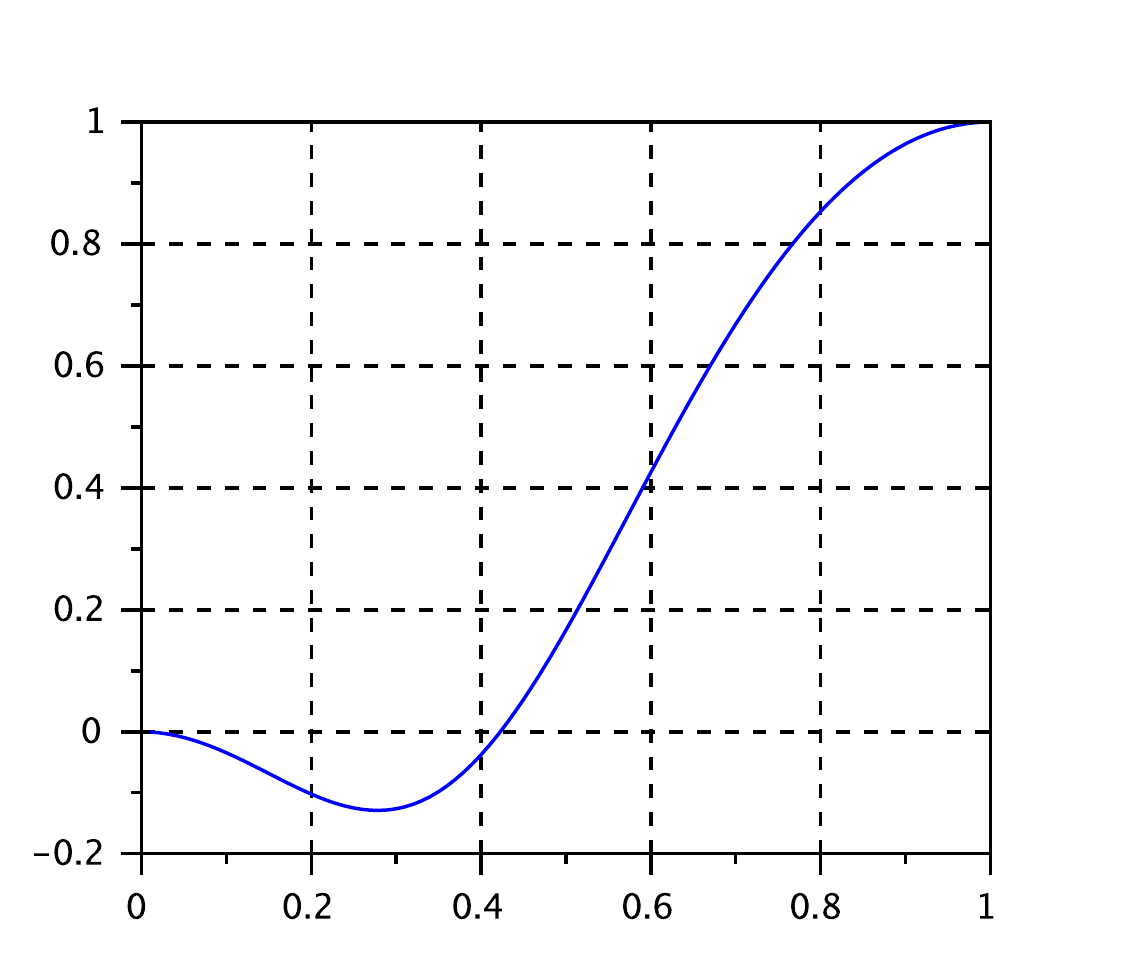}
\caption{A typical flux function $s\mapsto f(s,c,k)$ with gravitation force.}
\label{fig:f}
\end{center}
\end{figure}

We now consider the polymer flooding model \eqref{1.1}-\eqref{1.3}, 
taking into account the gravitation force but neglect the adsorption effect, 
i.e.,
\begin{align}
 &s_t + f(s,c,k)_x =0, \label{5.1} \\
& (cs)_t + (cf(s,c,k))_x =0,\label{5.2} \\
 &k_t =0.\label{5.3} 
\end{align}
An example for the flux function $f(s,c,k)$  was derived in~\cite{GimseRisebro3},  
where the flux function $f(s,c,k)$ typically becomes negative for small values of $s$, 
see Figure~\ref{fig:f}.
For simplicity of the discussion, we assume the monotone properties~\eqref{fconds}.

%%%%%%%%%%%%%%%%%%%%%%
\subsection{Lagrangian coordinates}

When we use the Lagrangian coordinate \eqref{eq:c},
it leads to the same system
\eqref{1.6}-\eqref{1.8}.
Let $A$ be the Jacobian matrix of the coordinate change, we have
$$
A=\begin{pmatrix} f & -s \\ 0 & 1
\end{pmatrix}, \qquad 
\det(A)=f.
$$
As $f$ changes the sign, $\det(A)$ changes sign as well, reversing the direction
of the ``time" variable in the Lagrangian system. 
Since such nonlinear PDE
is not time reversible, this coordinate change is not valid. 

In this case, 
one may introduce a modified Lagrangian coordinate $(\tilde\phi,\tilde\psi)$ as
\begin{equation}\label{Lag2}
\tilde \phi_x = - \mbox{sign}(f) \cdot s, \quad \tilde\phi_t =\mbox{sign}(f)\cdot f, 
\quad \tilde\phi(0,0)=0,
\qquad \tilde\psi=x.
\end{equation}
This leads to the following Lagrangian system:
\begin{equation}\label{5.6-5.8}
\left\{
\begin{array}{l} 
\displaystyle \left(\frac{1}{f}\right)_{\tilde\psi}
- \mbox{sign}(f) \cdot \left(\frac{s}{f}\right)_{\tilde\phi} =0,  \\[3mm]
\displaystyle c_{\tilde\psi} =0, \\[1mm]
\displaystyle k_{\tilde\phi}=0. 
\end{array}\right.
\end{equation}

%%%%%%%%%%%%%%%%%%%%%%%%%%%
\subsection{The reduced models}\label{sec5.2}
%%%%%%%%%%%%%%%%%%%%
There are two types of reduced models, for $k$=constant and for $c$=constant. 

\textbf{Type 1}.
When $k$ is constant, we have the reduced system~\eqref{1.10}-\eqref{1.20},
i.e.,
\begin{equation}\label{m1}
 \begin{cases} s_t + f(s,c)_x=0,\\
(cs)_t + (cf(s,c))_x=0,
\end{cases}
\end{equation}
where $s\mapsto f$ is as illustrated in Figure~\ref{fig:f}. 
The solution of the Riemann problem follows the same Riemann solver as in
section~\ref{sec2.1}, now with a different flux function $f(s,c)$.
We remark that, this reduced model was studied in detail in \cite{WS}, 
for a more general class of flux function $f(s,c)$, 
where existence of entropy solutions for the Cauchy problem was established. 

The solutions of the Riemann problem for~\eqref{m1} have the following properties.
Let $(s_l,c_l)$, $(s_r,c_r)$ be the Riemann data, and denote
$f_l=f(s_l,c_l)$. Let
$s_0>0 $ be the unique value such that $f(s_0,c_l)=0$.  
The followings hold.
\begin{itemize}
\item If $s_l < s_0$, i.e., $f_l<0$, 
then the $c$-wave in the solution travels with negative speed.  
\item If $s_l>s_0$, i.e., $f_l>0$, 
then the $c$-wave in the solution travels with positive speed.  
\item If $s_l =s_0$, i.e. $f_l=0$, then the $c$-wave is stationary. 
\end{itemize}
We remark that these properties give us the information on the ordering 
between the $c$  wave and the $k$ wave
in the Riemann solution of the $3\times3$ system,
making that Riemann solver in Section~\ref{sec5.3} 
easier to construct.

\medskip

\textbf{Type 2}.
When $c\equiv$ constant, we have the following $2\times2$ system
\begin{equation}\label{m2}
 \begin{cases} s_t + f(s,k)_x =0,\\
  k_t =0.
\end{cases}
\end{equation}
Since $f_s$ can change sign, the system is parabolic degenerate at $f_s=0$. 
The Riemann solver follows the same construction as for a scalar conservation
law with discontinuous flux, where the key step is to locate
the path of $k$ wave. 
Let $ f_l(s)=f(s,k_l)$ and $f_r(s)=f(s,k_r)$, and define $G^\sharp,G^\flat$ 
accordingly as in~\eqref{Gsharp}-\eqref{Gflat}.
The {\em minimum jump path} connecting $G^\sharp$ and $G^\flat$ 
is the $k$ wave. The rest follows.

%%%%%%%%%%%%%%%%%%%%%%%%%%%
\subsection{The Riemann solver for the $3\times 3$ Eulerian system}
\label{sec5.3}

Regardless of the signs of the wave speeds, we have the following properties:
\begin{itemize}
\item The $k$-wave travels with speed 0. Crossing it, $c$ and $f$ remain continuous;
\item Crossing a $c$-wave, $k$ and $f/s$ remain continuous;
\item Crossing an $s$-wave, $k$ and $c$ remain continuous.
\end{itemize}

The properties in Section~\ref{sec5.2} give us the sign of the speed for the  $c$ wave,
even for the full system. 
Let $(s_l,c_l,k_l)$, $(s_r,c_r,k_r)$ be the Riemann data, and let 
$f_l=f(s_l,c_l,k_l)$. 
Then, if $f_l<0$ or $s_l=0$,  the $c$-wave speed is negative;
 if $f_l>0$, the $c$ wave speed is positive. 
We can now construct  the Riemann solver. 

\textbf{Case (1): The $c$ wave travels with negative speed.} 
Let $(s_m, c_r,k_l)$ denote the trace at $x=0-$ in the solution of the Riemann problem.
We need to solve two Riemann problems: 
\begin{itemize}
\item[(R1):] Riemann problem between the states $(s_l,c_l,k_l)$ and $(s_m,c_r,k_l)$, i.e.,
$$ 
\begin{cases} s_t + f(s,c,k_l)_x=0, \\ (cs)_t + (cf(s,c,k_l))_x=0, 
\end{cases} \qquad 
(s,c)(0,x)=\begin{cases} (s_l,c_l), & (x<0) ,\\  (s_m,c_r), & (x>0) .\end{cases}
$$
This is a reduced model of type 1, discussed in Section~\ref{sec5.2}.
\item[(R2):] Riemann problem between the states $(s_m,c_r,k_l)$ and $(s_r,c_r,k_r)$, i.e.,
$$ 
\begin{cases} s_t + f(s,c_r,k)_x=0, \\ k_t=0, 
\end{cases}
\qquad 
(s,k)(0,x)=\begin{cases} (s_m,k_l), & (x<0) ,\\  (s_r,k_r), & (x>0). \end{cases}
$$
This is a reduced model of type 2, discussed in Section~\ref{sec5.2}.
\end{itemize}

For the solution to be plausible,  the speeds for the waves of (R1) must be $<0$, and
the speeds of the waves  of (R2) must be $\ge 0$. 
Here we use the strict ``$<$'' relation for the waves from (R1), 
to ensure that $s_m $ is the trace at $x=0-$,
rather than a middle state of two stationary waves. 

We denote various related flux functions as:
$$ 
f_l(s)\; \dot=\; f(s,c_l,k_l),  \qquad f_m(s) \;\dot=\; f(s,c_r,k_l), \qquad
 f_r(s) \;\dot=\; f(s,c_r,k_r).
$$
Given the value $s_l$ and the two flux functions $f_l,f_m$, 
let $I_1$ denote the set of values for $s_m$ such that the Riemann problem 
(R1) is solved with  waves of negative speed.  
Since the Riemann solution for the above system is 
uniquely defined, the set $I_1$ can be uniquely constructed,
following this general algorithm. 
For any given $s_l$, we locate all possible 
$c$~wave paths that  travel with negative speed.
Then, for all the $c$~wave paths, we locate all possible $s_m$ 
that connects to the $c$~wave with $s$~waves of negative speed.

\begin{figure}[htbp]
\begin{center}
\setlength{\unitlength}{1mm}
\begin{picture}(80,50)(0,10)   % -- top left
\put(-3,48){\small $ (f_l>f_m, s_l \ge \hat s)$}
\put(0,40){\line(1,0){60}}
\multiput(20,24)(0,1){16}{\line(0,1){0.5}}\put(19.2,41){\small $\hat s$}
\multiput(0,40)(1,-0.8){20}{\line(5,-4){0.5}}
\multiput(0,40)(1,-0.6){35}{\line(2,-1){0.5}}
\multiput(28,23.2)(0,1){17}{\line(0,1){0.5}}
\multiput(17,20.8)(1,0){16}{\line(1,0){0.5}}
\multiput(32.2,20.8)(0,1){19}{\line(0,1){0.5}}
\multiput(17.3,20.8)(0,1){19}{\line(0,1){0.5}}
\thicklines
\color{red}
\qbezier(0,40)(10,40)(15,30)\qbezier(15,30)(26,10)(52,45)
\put(48,45){\small $f_l$}
\put(27,41){\small $s_l$}
\put(28,23.2){\circle*{0.8}} 
\color{blue}
\qbezier(0,40)(8,40)(12,30)
\qbezier(12,30)(23,0)(55,45)\put(55,43){\small $f_m$}
\color{green} 
\qbezier(0,40)(8,40)(12,30)\qbezier(12,30)(14,24)(17.3,20.8)
\put(32.2,20.8){\circle*{0.8}}
\put(32,41){\small$s_2$}
\put(15,41){\small$s_1$}
\put(17.3,20.8){\circle*{0.8}}
\put(5,12){$I_1=(0,s_1)\cup\{s_2\}$}
\color{black}
\put(0,40){\circle*{0.8}}\put(20,24){\circle*{0.8}}
\end{picture}
\begin{picture}(60,50)(0,10)   % -- top right
\put(-3,48){\small $ (f_l>f_m, s_l \le \hat s)$}
\put(0,40){\line(1,0){60}}
\multiput(20,24)(0,1){16}{\line(0,1){0.5}}\put(18,37){\small $\hat s$}
\multiput(0,40)(1,-0.8){29}{\line(5,-4){0.5}}
\multiput(12,34.5)(0,1){6}{\line(0,1){0.5}}
\multiput(17,18.5)(1,0){16}{\line(1,0){0.5}}
\multiput(27,18.5)(0,1){22}{\line(0,1){0.5}}
\multiput(21,18.5)(0,1){22}{\line(0,1){0.5}}
\thicklines
\color{red}
\qbezier(0,40)(10,40)(15,30)\qbezier(15,30)(26,10)(52,45)\put(48,45){\small $f_l$}
\put(12,34.5){\circle*{0.8}} 
\put(10,41){\small $s_l$}
\color{blue}
\qbezier(0,40)(8,40)(12,30)
\qbezier(12,30)(23,0)(55,45)\put(55,43){\small $f_m$}
\color{green} 
\qbezier(0,40)(8,40)(12,30)\qbezier(12,30)(15.2,20.8)(21,18.5)
\put(27,18.5){\circle*{0.8}}
\put(26,41){\small$s_2$}
\put(21,41){\small$s_1$}
\put(21,18.5){\circle*{0.8}} 
\put(5,12){$I_1=(0,s_1)\cup\{s_2\}$}
\color{black}
\put(0,40){\circle*{0.8}}
\put(20,24){\circle*{0.8}}
\end{picture}
\begin{picture}(80,50)(0,10)   % -- bottom left
\put(-3,48){\small $ (f_l>f_m, s_l \ge \tilde s)$}
\put(0,40){\line(1,0){60}}
\multiput(25,22.6)(0,1){18}{\line(0,1){0.5}}\put(24,41){\small $\hat s$}
\multiput(0,40)(1,-0.69){31}{\line(5,-4){0.5}}
\multiput(0,40)(1,-0.4){40}{\line(2,-1){0.5}}
\multiput(37.8,24.8)(0,1){16}{\line(0,1){0.5}}
\multiput(29.7,19.4)(0,1){21}{\line(0,1){0.5}}\put(29,41){\small $\tilde s$}
\multiput(18,26.2)(1,0){17}{\line(1,0){0.5}}
\multiput(34.2,26.2)(0,1){14}{\line(0,1){0.5}}
\multiput(17.6,26.2)(0,1){14}{\line(0,1){0.5}}
\thicklines
\color{blue}
\qbezier(0,40)(10,40)(15,30)\qbezier(15,30)(26,10)(52,45)\put(47,45){\small $f_m$}
\color{red}
\qbezier(0,40)(8,40)(12,30)
\qbezier(12,30)(23,0)(55,45)\put(55,43){\small $f_l$}
\put(37.8,24.8){\circle*{0.8}} 
\put(37,41){\small $s_l$}
\color{green} 
\qbezier(0,40)(10,40)(15,30)
\qbezier(15,30)(16,28)(17.6,26.2)
\put(34.2,26.2){\circle*{0.8}}
\put(15,41){\small$s_1$}
\put(33,41){\small$s_2$}
\put(17.6,26.2){\circle*{0.8}}
\put(5,12){$I_1=(0,s_1)\cup\{s_2\}$}
\color{black}
\put(0,40){\circle*{0.8}}
\put(25,22.6){\circle*{0.8}} 
\put(29.7,19.4){\circle*{0.8}}
\end{picture}
\begin{picture}(60,50)(0,10)   % -- bottom right
\put(-3,48){\small $ (f_l>f_m, s_l \le \tilde s)$}
\put(0,40){\line(1,0){60}}
\multiput(25,22.6)(0,1){18}{\line(0,1){0.5}}\put(24,41){\small $\hat s$}
\multiput(0,40)(1,-0.69){31}{\line(5,-4){0.5}}
\multiput(29.7,19.4)(0,1){21}{\line(0,1){0.5}}\put(29,41){\small $\tilde s$}
\multiput(14.5,24.8)(0,1){16}{\line(0,1){0.5}}
\thicklines
\color{blue}
\qbezier(0,40)(10,40)(15,30)\qbezier(15,30)(26,10)(52,45)\put(47,45){\small $f_m$}
\color{red}
\qbezier(0,40)(8,40)(12,30)
\qbezier(12,30)(23,0)(55,45)\put(55,43){\small $f_l$}
\put(14.5,24.8){\circle*{0.8}}
\put(14,41){\small $s_l$}
\color{green}
\qbezier(0,40)(10,40)(15,30)
\qbezier(15,30)(18.5,23)(25,22.6)
\put(5,12){$I_1=(0,\hat s]$}
\color{black}
\put(0,40){\circle*{0.8}}
\put(25,22.6){\circle*{0.8}} 
\put(29.7,19.4){\circle*{0.8}}
\end{picture}
\caption{The set $I_1$ for various  cases. }
\label{fig:I1}
\end{center}
\end{figure}
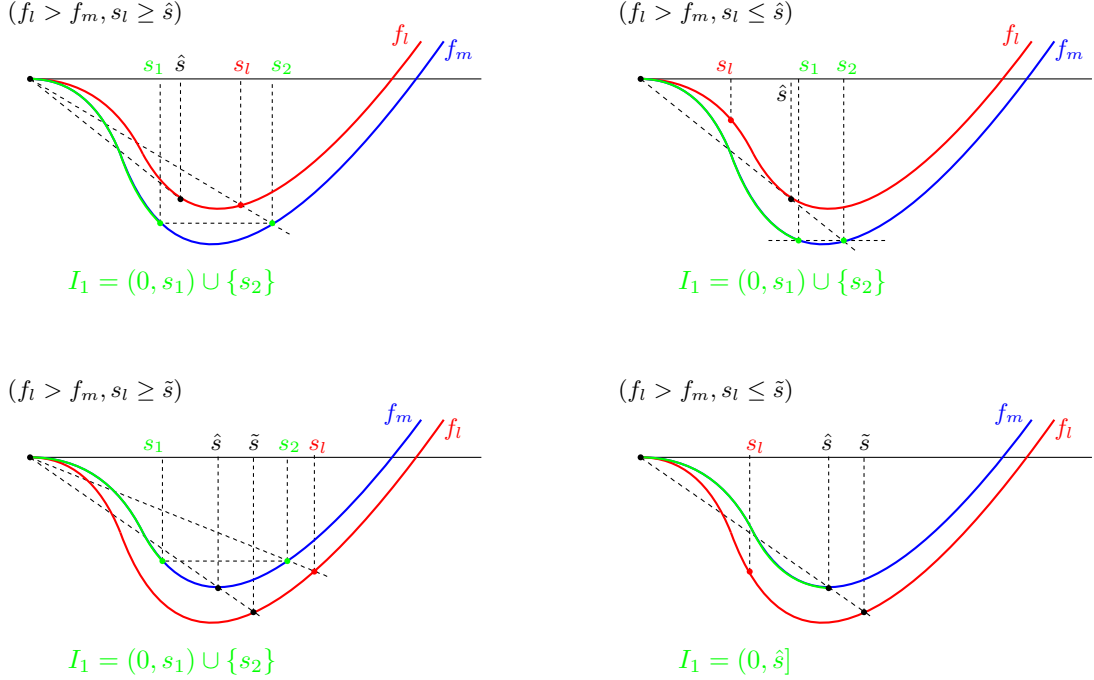

There are 4 cases, illustrated in Figure~\ref{fig:I1}.
\begin{itemize}
\item Consider $c_l<c_r$, and therefore $f_l>f_m$. 
Let $\hat s>0$ be the unique value that satisfies 
$f_l(\hat s)/\hat s=  f_l'(\hat s)$.
The set $I_1$ depends on the relation between $s_l$ and $\hat s$. 
\begin{itemize}
\item 
Consider $s_l \ge \hat s$.  We discuss only this case in detail, as an explanation  of the
general algorithm, while the other cases being similar.  
In Figure~\ref{fig:I1detail}, 
the location of $s_l$ is marked in red. 
The value $s_2>s_l$ is chosen such that the three points
$(0,0), (s_l,f_l(s_l)),(s_2,f_m(s_2))$ are collinear.
The points $s_3,s_4$ lie on the same line. 
We also have $s_1\le s_2$ with
$ f_m(s_1)=f_m(s_2)$. 
There are two sub-cases.

(1) On the right side of $\hat s$, the $c$ wave path can only be the one connecting 
$s_l$ and $s_2$. Thus, $\{s_2\} \in I_1$. 
To connect it further with negative $s$ waves, one can connect $s_2$ to 
any $s$ value between $s_4$ and $s_1$ with a $s$ shock. 
Thus $[s_4,s_1)\in I_1$.  

(2) Now we consider the case where the $c$ wave is on the left side of $\hat s$. 
Then, $s_l$ can be connected to an $s<s_3$ with an $s$ shock,
then connects to $f_m$ through a $c$ wave, reaching a point on the left of $s_4$.
This point can further be connected to even smaller values of $s$ through 
$s$ waves.  Thus $(0, s_4)\in I_1$.

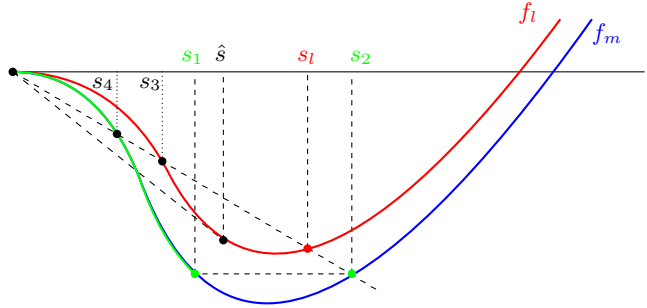
\begin{figure}[htbp]
\begin{center}
\setlength{\unitlength}{1.4mm}
\begin{picture}(50,33)(0,18)   % -- top left
\put(0,40){\line(1,0){60}}
\multiput(20,24)(0,1){16}{\line(0,1){0.5}}\put(19.2,41){\small $\hat s$}
\multiput(0,40)(1,-0.8){20}{\line(5,-4){0.5}}
\multiput(0,40)(1,-0.6){35}{\line(2,-1){0.5}}
\multiput(28,23.2)(0,1){17}{\line(0,1){0.5}}
\multiput(17,20.8)(1,0){16}{\line(1,0){0.5}}
\multiput(32.2,20.8)(0,1){19}{\line(0,1){0.5}}
\multiput(17.3,20.8)(0,1){19}{\line(0,1){0.5}}
\thicklines
\color{red}
\qbezier(0,40)(10,40)(15,30)\qbezier(15,30)(26,10)(52,45)
\put(48,45){\small $f_l$}
\put(27,41){\small $s_l$}
\put(28,23.2){\circle*{0.8}} 
\color{blue}
\qbezier(0,40)(8,40)(12,30)
\qbezier(12,30)(23,0)(55,45)\put(55,43){\small $f_m$}
\color{green} 
\qbezier(0,40)(8,40)(12,30)\qbezier(12,30)(14,24)(17.3,20.8)
\put(32.2,20.8){\circle*{0.8}}
\put(32,41){\small$s_2$}
\put(16,41){\small$s_1$}
\put(17.3,20.8){\circle*{0.8}}
\color{black}
\put(0,40){\circle*{0.8}}\put(20,24){\circle*{0.8}}
\put(14.2,31.5){\circle*{0.8}}
\put(9.9,34.1){\circle*{0.8}}
\thinlines
\qbezier[20](14.2,31.5)(14.2,35.7)(14.2,40)
\qbezier[17](9.9,34.1)(9.9,37)(9.9,40)
\put(12,38.5){\small $s_3$}\put(7.5,38.5){\small $s_4$}
\end{picture}
\caption{The set $I_1$ for  case 1, with details for the construction. }
\label{fig:I1detail}
\end{center}
\end{figure}

This case is summarized in the top-left plot in Figure~\ref{fig:I1}. 
The set $I_1$ contains an open interval $(0,s_1)$ and an isolated point $s_2$, 
i.e.,
$$ I_1 = (0,s_1)\cup \{s_2\}.$$ 

\item 
Consider $s_l \le \hat s$, and see the top-right plot in Figure~\ref{fig:I1}. 
Here,  the three points $(0,0), (\hat s, f_l(\hat s)), (s_2, f_m(s_2))$ are collinear, 
and we have $s_1<s_2$ and $f_m(s_1)=f_m(s_2)$. 
This gives 
$  I_1 = (0,s_1)\cup \{s_2\}$. 
\end{itemize}
\item
Consider $c_l>c_r$ and therefore $f_l < f_m$.
Let $\hat s >0$ be the point where $f_m'(\hat s)=0$, and 
let $\tilde s>0$ be the unique point such that 
the three points $(0,0), (\hat s, f_m(\hat s)), (\tilde s, f_l(\tilde s))$ are collinear. 
The set $I_1$ depends on the relation between $s_l$ and $\tilde s$. 
We have 2 cases.
\begin{itemize}
\item
Consider $s_l \ge \tilde s$, and see the bottom-left plot in Figure~\ref{fig:I1}. 
Here $s_2>0$ is the chosen such that the three points 
$(0,0), (s_l, f_l(s_l)), (s_2,f_m(s_2))$ are collinear. 
The value $s_1$ is chosen such that $s_1<s_2$ and 
$f_m(s_1)=f_m(s_2)$. 
  Then we have 
$$  I_1 = (0,s_1)\cup \{s_2\}.$$ 
\item
Consider $s_l \le \tilde s$, and see the bottom-right plot in Figure~\ref{fig:I1}. 
For this case we simply have 
$ I_1=(0,\hat s]$. Note that $\hat s\in I_1$ since any $s$-wave connecting to $\hat s$
must be a rarefaction. 
\end{itemize}
\end{itemize}

\medskip

\begin{figure}[htbp]
\begin{center}
\setlength{\unitlength}{1mm}
\begin{picture}(80,50)(0,-30)   % -- top left
\put(-3,10){\small $ (f_m>f_r, s_r \ge \tilde s)$}
\put(0,0){\line(1,0){60}}
\multiput(8,-15)(1,0){16}{\line(1,0){0.5}}
\multiput(25,-17)(0,1){17}{\line(0,1){0.5}}
\multiput(12,-15)(0,1){15}{\line(0,1){0.5}}\put(11,1){\small $\tilde s$}
\multiput(20.2,-15)(0,1){15}{\line(0,1){0.5}}\put(20,1){\small $\hat s$}
\thicklines
\color{blue}
\qbezier(0,0)(10,0)(15,-10)
\qbezier(15,-10)(20,-20)(26,-10)
\qbezier(26,-10)(40,15)(55,15)
\put(48,15.5){\small $f_m$}
\color{red}
\qbezier(0,0)(8,0)(12,-15)
\qbezier(12,-15)(17,-32)(26,-15)
\qbezier(26,-15)(40,15)(55,15)
\put(53,12){\small $f_r$}
\put(25,-17){\circle*{0.8}} 
\put(24,1){\small $s_r$}
\color{green} 
\qbezier(26,-10)(40,15)(55,15)
\qbezier(20,-15)(23.2,-15)(26,-10)
\put(35,-20){$I_2=[\hat s,1]$}
\color{black}
\put(12,-15){\circle*{0.8}} 
\put(20.2,-15){\circle*{0.8}} 
\end{picture}
\begin{picture}(60,50)(0,-30)   % -- top right
\put(-3,10){\small $ (f_m>f_r, s_r \le \tilde s)$}
\put(0,0){\line(1,0){60}}
\multiput(8,-15)(1,0){16}{\line(1,0){0.5}}
\multiput(9.5,-8)(0,1){8}{\line(0,1){0.5}} 
\multiput(9.5,-8)(1,0){18}{\line(1,0){0.5}} 
\multiput(12,-15)(0,1){15}{\line(0,1){0.5}}\put(11,1){\small $\tilde s$}
\multiput(20.2,-15)(0,1){15}{\line(0,1){0.5}}\put(20,1){\small $\hat s$}
\multiput(13.9,-8)(0,1){8}{\line(0,1){0.5}} 
\multiput(27.1,-8)(0,1){8}{\line(0,1){0.5}} 
\thicklines
\color{blue}
\qbezier(0,0)(10,0)(15,-10)
\qbezier(15,-10)(20,-20)(26,-10)
\qbezier(26,-10)(40,15)(55,15)
\put(48,15.5){\small $f_m$}
\color{red}
\qbezier(0,0)(8,0)(12,-15)
\qbezier(12,-15)(17,-32)(26,-15)
\qbezier(26,-15)(40,15)(55,15)
\put(53,12){\small $f_r$}
\put(9.5,-8){\circle*{0.8}} 
\put(7,1){\small $s_r$}
\color{green}
\put(13.9,-8){\circle*{0.8}} \put(13,1){\small $s_1$} 
\put(27.1,-8){\circle*{0.8}} \put(26,1){\small $s_2$} 
\qbezier(27.1,-8)(40.5,15)(55,15)
\put(35,-20){$I_2=[s_2,1]\cup \{s_1\}$}
\color{black}
\put(12,-15){\circle*{0.8}}
\put(20.2,-15){\circle*{0.8}}
\end{picture}
\begin{picture}(80,50)(0,-30)   % -- bottom left
\put(-3,10){\small $ (f_m<f_r, s_r \ge \hat s)$}
\put(0,0){\line(1,0){60}}
\multiput(8,-15)(1,0){19}{\line(1,0){0.5}}
\multiput(12,-15)(0,1){15}{\line(0,1){0.5}} 
\multiput(26,-15)(0,1){15}{\line(0,1){0.5}} 
\multiput(20.2,-15)(0,1){15}{\line(0,1){0.5}}\put(20,1){\small $\hat s$}
\thicklines
\color{red}
\qbezier(0,0)(10,0)(15,-10)
\qbezier(15,-10)(20,-20)(26,-10)
\qbezier(26,-10)(40,15)(55,15)
\put(48,15.5){\small $f_r$}
\put(34.5,3){\circle*{0.8}} 
\put(31,3){\small $s_r$}
\color{blue}
\qbezier(0,0)(8,0)(12,-15)
\qbezier(12,-15)(17,-32)(26,-15)
\qbezier(26,-15)(40,15)(55,15)
\put(53,12){\small $f_m$}
\color{green} 
\qbezier(26,-15)(40,15)(55,15)
\put(12,-15){\circle*{0.8}} \put(11,1){\small $s_1$} 
\put(26,-15){\circle*{0.8}} \put(25,1){\small $s_2$} 
\put(35,-20){$I_2=[s_2,1]\cup\{s_1\}$}
\color{black}
\put(20.2,-15){\circle*{0.8}} 
\end{picture}
\begin{picture}(60,50)(0,-30)   % -- bottom right
\put(-3,10){\small $ (f_m<f_r, s_r \le \hat s)$}
\put(0,0){\line(1,0){60}}
\multiput(8,-15)(1,0){19}{\line(1,0){0.5}}
\multiput(8,-8)(1,0){23}{\line(1,0){0.5}}
\multiput(9.6,-8)(0,1){8}{\line(0,1){0.5}} 
\multiput(29.5,-8)(0,1){8}{\line(0,1){0.5}} 
\multiput(13.9,-8)(0,1){8}{\line(0,1){0.5}}
\multiput(20.2,-15)(0,1){15}{\line(0,1){0.5}}\put(20,1){\small $\hat s$}
\thicklines
\color{red}
\qbezier(0,0)(10,0)(15,-10)
\qbezier(15,-10)(20,-20)(26,-10)
\qbezier(26,-10)(40,15)(55,15)
\put(48,15.5){\small $f_r$}
\put(13.9,-8){\circle*{0.8}} \put(13,1){\small $s_r$} 
\color{blue}
\qbezier(0,0)(8,0)(12,-15)
\qbezier(12,-15)(17,-32)(26,-15)
\qbezier(26,-15)(40,15)(55,15)
\put(53,12){\small $f_m$}
\color{green} 
\qbezier(29.5,-8)(41.8,15)(55,15)
\put(9.6,-8){\circle*{0.8}} \put(8.5,1){\small $s_1$} 
\put(29.5,-8){\circle*{0.8}} \put(28.5,1){\small $s_2$} 
\put(35,-20){$I_2=[s_2,1]\cup\{s_1\}$}
\color{black}
\put(20.2,-15){\circle*{0.8}} 
\end{picture}
\caption{Illustrations for the set $I_2$ for various  cases.} 
\label{fig:I2}
\end{center}
\end{figure}
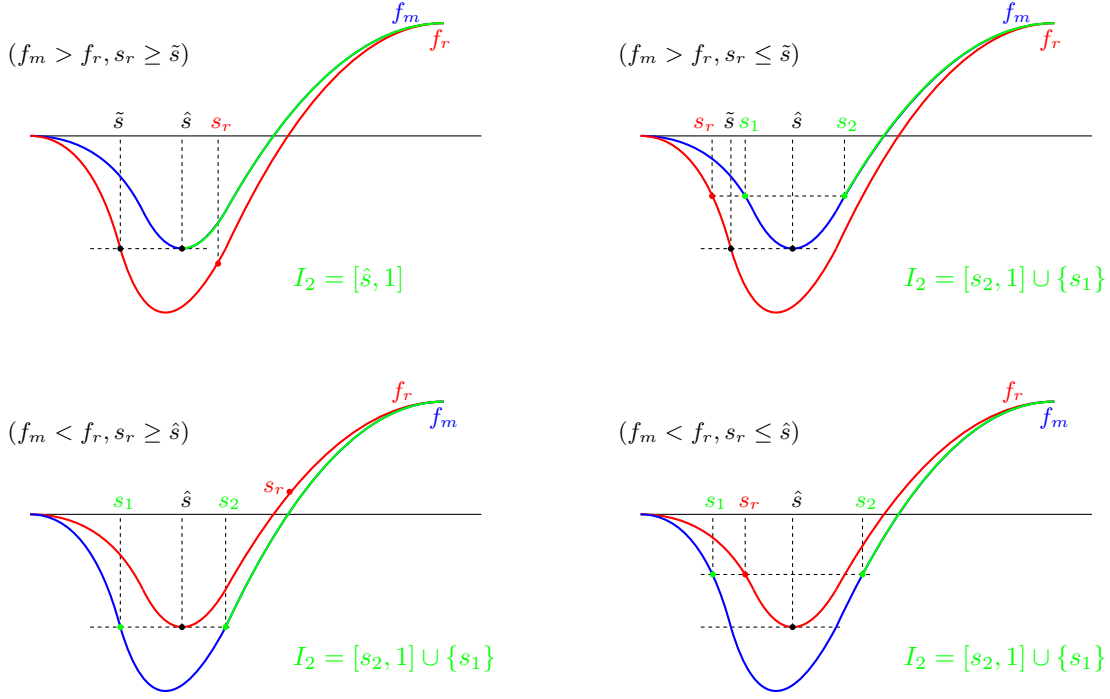

Now we consider problem (R2), for $x\ge 0$. 
Given the value $s_r$ and the flux functions $f_m,f_r$, 
we denote by $I_2$ the set of values for $s_m$ such that (R2)
is solved with non-negative speed. 
There are 4 cases, illustrated in Figure~\ref{fig:I2}. 
\begin{itemize}
\item If $k_l>k_r$, we have $f_m>f_r$. 
Let $\hat s$ be the resonant point where $f_m'(\hat s)=0$,
and $\tilde s<\hat s$ satisfies $f_r(\tilde s)=f_m(\hat s)$. 
We have 2 sub-cases.
  \begin{itemize}
  \item
  For $s_r\ge\tilde s$, we illustrate it in the top-left plot in Figure~\ref{fig:I2}. 
  Here we simply have 
  $$I_2 = [\hat s, 1].$$ 
  \item
  For $s_r\le \tilde s$, we illustrate it in the top-right plot in Figure~\ref{fig:I2}. 
  Here $s_1<s_2$ are two unique values such that 
  $ f_m(s_1)=f_m(s_2)=f_r(s_r)$.  The set $I_2$ contains a closed set $[s_2,1]$
  and an isolated point $s_1$, i.e.,
  $$ I_2=[s_2,1]\cup \{s_1\}.$$
  \end{itemize}

\item If $k_l<k_r$, we have $f_m < f_r$. 
Let $\hat s>0$ be the resonant point such that $f_r'(\hat s)=0$. 
Again we have 2 sub-cases.
  \begin{itemize}
  \item The case $s_r \ge \hat s$ is illustrated in the bottom-left plot of Figure~\ref{fig:I2}.
  Here $s_1<s_2$ are chosen such that 
  $ f_m(s_1)=f_r(\hat s)=f_m(s_2)$.  We get 
  $$ I_2=[s_2,1]\cup \{s_1\}.$$
  \item
  The case $s_r \le \hat s$ is illustrated in the bottom-right plot of Figure~\ref{fig:I2}.
  Here $s_1<s_2$ are chosen such that 
  $ f_m(s_1)=f_r(s_r)=f_m(s_2)$. 
  We get 
  $$ I_2=[s_2,1]\cup \{s_1\}.$$
  \end{itemize}
\end{itemize}

\bigskip

We now summarize. 
 For all cases, the flux $f_m$ is decreasing on the set $I_1$
and increasing on the set $I_2$. 
  Furthermore, there exits only one single point
on the set $I_1$ where $(f_m)' \ge 0$, 
and one single point on the set $I_2$ 
where $(f_m)' \le 0$. 
Thus, for any combination of the pair $I_1,I_2$, the intersection 
$I_1 \cap I_2$ is non-empty and consists of a exactly one single point.  
We now let the $s$ value of  this point be  the trace $s_m$.  
Note that in all cases, we have $f_m(s_m)<0$. 

\bigskip

\textbf{Case (2):  The $c$ wave travels with positive speed.} 
Let $(s_m,c_l,k_r)$ be the trace along $x=0+$. 
We have two Riemann problems:
\begin{itemize}
\item[(R3):] Riemann problem connecting states $(s_l,c_l,k_l)$ and $(s_m,c_l,k_r)$,
which is a reduced model of type 2, and should be 
solved with waves of speed $\le 0$.
\item[(R4):] Riemann problem connecting states  $(s_m,c_l,k_r)$ and $(s_r,c_r,k_r)$,
which is a reduced model of type 1, and should be 
solved with waves of speed $> 0$.
\end{itemize}

With some abuse of notations, we denote the flux functions
$$
f_l(s)\;\dot=\; f(s, c_l,k_l), \qquad f_m(s)\;\dot=\; f(s, c_l,k_r).
$$
Let $I_3$ be the set for the $s_m$ values such that (R3) is solved with waves of
non-positive speed.
Recall that $f_l(s_l)>0$. Note then,
if $f_l>0$ and $f_m>0$, then $f'_l>0$ and $f'_m>0$. 
Thus, $I_3$ consists of a single point,  call it $s_m$, 
such that $ f_l(s_l)=f_m(s_m)>0$. 
As a result, the solution to (R3) consists  of a single stationary $k$-wave. 
 
It can be easily verified that with this $s_m$, Riemann problem
(R4) is solved with waves of positive speed.

\medskip

\textbf{Case (3):  The $c$ wave is stationary.} 
We have $f(s_l,c_l,k_l)=0$ and $s_l>0$. 
Let $ s_m>0$ be the unique value such that $f( s_m,c_r,k_r)=0$. 
We have a combined $c$+$k$ stationary wave 
connecting $(s_l,c_l,k_l)$ to $( s_m,c_r,k_r)$. 
Then, $( s_m,c_r,k_r)$ can be connected $(s_r,c_r,k_r)$ by solving
a Riemann problem with flux $f_r(s)=f(s,c_r,k_r)$. 
Thanks to the special location of $ s_m$, 
the solution consists of waves of non-negative speeds.

\medskip

In summary, we have a global Riemann solver which 
generates a unique solution for any initial Riemann data. 
Furthermore, all discontinuities are entropy admissible, i.e., they are limits
of vanishing viscosity solutions of viscous travel waves.

%%%%%%%%%%%%%%%%%%%%%%%%%%%%
\section{Concluding Remarks}\label{sec6}
\setcounter{equation}{0}

In this paper we construct global Riemann solvers for several $3\times3$ 
systems of conservations laws, arising in polymer flooding and traffic flow. 
However, we neglected the polymer flooding model where 
both the gravitation force and the adsorptive effect are considered.
For system~\eqref{0.1}-\eqref{0.3}
where $s\mapsto f$ is as in Figure~\ref{fig:f}, the analysis is more complicated,
since the Riemann solver for the reduced systems of Type 1 is not 
available in the literature.
Nevertheless, a global Riemann solver can be constructed, following a 
somewhat similar approach. Due to the new details involved, 
it deserves to be treated in a separate paper in near future.

It is also interesting to consider the case of multi-component polymer flooding.
For the system~\eqref{0.1}-\eqref{0.3}, 
the equation~\eqref{0.2} is replaced by an $n\times n$ system, where $n$
is the number of different types of polymers.  
The size of the full system is $(n+2)\times(n+2)$. 
We denote the families as the $\{s, c_1, \cdots, c_n, k\}$ families. 
For the non-adsorptive case where $m(c)=\mbox{constant}$, 
all the $c_i$ families are linearly degenerate, where all waves 
travel with non-negative speed and they never interact. 
A global Riemann solver can be easily constructed in a similar way as 
the one in Section~\ref{sec2} or Section~\ref{sec4}. 
For the adsorptive model without gravitation force, it depends heavily
on the adsorptive function $m(c)$, where $c \in R^n$ and the function $m(c)$ is
vector-valued. 
In the literature, various adsorptive functions have been studied, 
leading to very different systems of~\eqref{0.2}. 
Using the Langmuir isotherm (cf.~\cite{Rhee}), 
$$ m_i(c_1,c_2,\cdots, c_n) = 
\frac{\kappa_ic_i}{1+ \kappa_1c_1+\cdots+\kappa_nc_n},\qquad
i=1,2,\cdots,n,$$
the systems~\eqref{0.2} is a Temple class when $s$ is constant.
A rather simple construction for the global Riemann solver can still be achieved
following a similar algorithm as in Section~\ref{sec3} and 
utilizing the reduced Riemann solver
in~\cite{Dahl}.

\end{document}